\documentclass[10pt]{article}
\usepackage{latexsym}
\usepackage{fullpage}
\usepackage{amssymb}
\begin{document}
%
%
\def\be{\begin{equation}}
\def\ee{\end{equation}}
\def\bea{\begin{eqnarray}}
\def\eea{\end{eqnarray}}
%
%
\def\g{\gamma}
\def\G{\Gamma}
\def\d{\delta}
\def\t{\tau}
\def\s{\sigma}
\def\S{\Sigma}
\def\k{\kappa}
\def\l{\lambda}
\def\o{\omega}
\def\O{\Omega}
\def\z{\zeta}
\def\ra{\rangle}
\def\la{\langle}
\def\M{{\widehat M}}
%
%
\def\Z{\mathbb{Z}}
\def\R{\mathbb{R}}
\def\C{\mathbb{C}}
\def\H{\mathbb{H}}
\def\RP#1{\ifmmode\R P^{#1}\else$\R P^{#1}$\fi}
\def\CP#1{\ifmmode\C P^{#1}\else$\C P^{#1}$\fi}
\def\HP#1{\ifmmode\H P^{#1}\else$\H P^{#1}$\fi}
%
%
\def\SO{\mathbf{SO}}
\def\Sp{\mathbf{Sp}}
\def\Spin{\mathbf{Spin}}
\def\so{\mathfrak{so}}
\def\sp{\mathfrak{sp}}
\def\spin{\mathfrak{spin}}
\def\Cl{\mathbf{Cl}}
\def\CCl{\mathbb{C}\mathbf{l}}
\def\L{\Lambda}
\def\Sym{{\rm Sym}}
\def\End{{\rm End}}
\def\Aut{{\rm Aut}}
\def\Re{\mathrm{Re\;}}
\def\Im{\mathrm{Im\;}}
\def\id{{\rm id}\,}
\def\Ad{{\rm Ad}}
\def\ad{{\rm ad}}
%
%
\def\spb{\mathbf{S}}
\def\HB{\mathbf{H}}
\def\EB{\mathbf{E}}
\def\FB{\mathbf{F}}
\def\1{\mathbf{1}}
\def\2{\mathit{1}}
\def\I{\mathbf{I}}
\def\J{\mathbf{J}}
\def\K{\mathbf{K}}
\def\#{\sharp}
\def\b{\flat}
\def\es{\lrcorner\;}
\def\esn{\lrcorner_\circ}
\def\dn{\wedge_\circ}
\def\tfrac#1#2{{\textstyle\frac{#1}{#2}}}
\def\qed{\ifmmode\quad\Box\else$\quad\Box$\fi}
%
%
%
\def\proof{{\bf Proof.}\hskip2mm}
\def\definition{{\bf Definition.}\hskip2mm}
\def\example{{\bf Example.}\hskip2mm}
\def\leer{\vskip 2.7mm \noindent}
\renewcommand{\theequation}{\thesection.\arabic{equation}}
\newtheorem{Lemma}[equation]{Lemma}
\newtheorem{Proposition}[equation]{Proposition}
\newtheorem{Corollary}[equation]{Corollary}
\newtheorem{Remark}[equation]{Remark}
\newtheorem{Theorem}{Theorem}

\title{The First Eigenvalue of the Dirac
 Operator on Quaternionic K\"ahler Manifolds}
\author{
 W. Kramer\thanks{Supported by the SFB 256 `Nichtlineare partielle
  Differentialgleichungen'},
 U. Semmelmann\thanks{Supported by the Max--Planck--Institut f\"ur
  Mathematik}, 
 G. Weingart\thanks{Supported by the SFB 256 `Nichtlineare partielle
  Differentialgleichungen'}\\
 {\small Mathematisches Institut der Universit\"at Bonn}\\
 {\small Beringstra\ss{}e 1, 53115 Bonn, Germany}\\
 {\small\texttt{kramer@math.uni-bonn.de,
   uwe@math.uni-bonn.de, gw@math.uni-bonn.de}}}
\maketitle
\begin{abstract}
 In \cite{qklast} we proved a lower bound for the spectrum of the Dirac 
 operator on quaternionic K\"ahler manifolds. In the present article we
 show that the only manifolds in the limit case, i.~e.~the only
 manifolds where the lower bound is attained as an eigenvalue, are the
 quaternionic projective spaces. We use the equivalent formulation in
 terms of the quaternionic Killing equation introduced in \cite{qklim} and
 show that a nontrivial solution defines a parallel spinor on the associated
 hyperk\"ahler manifold.
\end{abstract}

\vskip0.5cm

\tableofcontents

\newpage

\section{Introduction}

 The square of the first eigenvalue of the Dirac operator on the sphere $S^n$
 of scalar curvature $\k$ is $\tfrac\k{4}\tfrac{n}{n-1}$. In \cite{fri5}
 Friedrich showed that this eigenvalue is a universal lower bound for all
 eigenvalues on an arbitrary compact Riemannian spin manifold $(M^n,\,g)$
 with positive scalar curvature $\k$ in the following sense:
 all eigenvalues of the Dirac operator satisfy
 $$
   \l^2 \ge {\mathrm{min}_M\k\over 4}\,{n\over n-1}\,.
 $$
 An eigenspinor realizing this lower bound is characterized by a special
 differential equation called the Killing equation. Conversely, on manifolds
 admitting Killing spinors, i.~e.~nontrivial solutions of this Killing 
 equation,  the lower bound is realized as an eigenvalue. These manifolds
 have been characterized by C.~B\"ar \cite{baer} translating the Killing
 equation on $M$ into the equation of a parallel spinor on the cone
 $\M=\R^+\times_{t^2}M$. The existence problem of solutions of
 the Killing equation was thus reduced to the description of manifolds
 with parallel spinors by M.~Wang \cite{wang}.

 Despite the fact that Friedrich's estimate is sharp, it is not optimal
 if $M$ is assumed to have additional geometric structure, namely special
 holonomy. Due to a result of O.~Hijazi and A.~Lichnerowicz
 (\cite{hijazi84}, \cite{lic}), there is no solution of the Killing
 equation if $M$ possesses a non--trivial parallel $k$--form, $k\neq 0, n$.
 There are two canonical classes of such manifolds which in addition have
 positive scalar curvature: K\"ahler manifolds and quaternionic K\"ahler
 manifolds. 

 The eigenvalue estimate for K\"ahler manifolds has been improved by
 K.--D.~Kirchberg in \cite{kir2} and \cite{kir6} (see also \cite{lic1}
 and \cite{hij1}). Again, this estimate is sharp: the lower bound is
 attained as first eigenvalue on the complex projective space $\CP{m}$
 resp.~on its product with a flat 2--torus in odd resp.~even complex
 dimensions. On K\"ahler manifolds of odd complex dimension, a spinor
 with smallest possible eigenvalue is characterized by a suitable
 modification of the Killing equation, the K\"ahlerian Killing equation.
 A.~Moroianu showed in \cite{andrei} that a K\"ahlerian Killing spinor
 defines an ordinary Killing spinor on the canonical $S^1$--bundle over
 $M$. Hence, the holonomy argument of B\"ar's work can be used to study
 the limit case. In even complex dimensions, the problem of characterizing
 the limit case is still open. Nevertheless, there are partial results
 by A.~Moroianu \cite{moroianu2} and A.~Lichnerowicz \cite{lic1}.

 A quaternionic K\"ahler manifold is an oriented $4n$--dimensional Riemannian
 manifold with $n\geq 2$ whose holonomy group is contained in the subgroup
 $\Sp(1)\cdot\Sp(n)\subset\SO(4n)$. Equivalently they are characterized by
 the existence of a certain parallel 4--form $\O$, the so--called fundamental
 or Kraines form (cf.~\cite{bon}, \cite{krain}). All quaternionic K\"ahler
 manifolds are Einstein with constant scalar curvature (cf.~\cite{alex2a}
 or \cite{ish3}) and possess a unique spin structure if $n$ is even,
 whereas for odd $n$ only the quaternionic projective spaces are spin
 (cf.~\cite{sala}). In \cite{qklast} we proved that on a compact quaternionic
 K\"ahler spin manifold $(M^{4n},\,g)$ of positive scalar curvature $\k$ the
 eigenvalues $\l$ of the Dirac operator satisfy
 $$
   \l^2 \ge {\k\over 4}\,{n+3\over n+2}\,.
 $$
 As in the Riemannian or K\"ahler case, this estimate is sharp, and the lower
 bound is attained as first eigenvalue on the quaternionic projective
 space (cf.~\cite{mil}).

 The natural task is to study the limit case and to find all manifolds
 which have $\frac{\k}{4}\,\frac{n+3}{n+2}$ in the spectrum of $D^2$. A
 first step in this direction was taken in \cite{qklim}, where we introduced
 the equation characterizing an eigenspinor with this particular eigenvalue.
 A new feature of this quaternionic Killing equation is that it involves not
 only the eigenspinor but also an auxiliary section of an additional bundle,
 which is not itself a spinor. We used it to show that no compact symmetric
 quaternionic K\"ahler manifolds besides the quaternionic projective spaces
 carry quaternionic Killing spinors. In the present article we prove the
 following more general
 \setcounter{Theorem}{0}
 \begin{Theorem}
 Let $M$ be a compact quaternionic K\"ahler manifold of quaternionic
 dimension $n$ and positive scalar curvature $\k$. If there is an
 eigenspinor for the Dirac operator with eigenvalue $\l$ satisfying
 $$
   \l^2 = \frac{\k}{4}\frac{n+3}{n+2}\,,
 $$
 then $M$ is isometric to the quaternionic projective space.
 \end{Theorem}
 For the proof we follow the approach of C.~B\"ar and A.~Moroianu.
 We consider the canonical $\SO(3)$--bundle $S$ associated with any
 quaternionic K\"ahler manifold. Introducing an appropriate metric on the
 total space $S$ the
 warped product $\M:=\R^+\times_{t^2}S$ has a natural hyperk\"ahler
 metric. Reformulating the Killing equation in terms of equivariant functions
 on the $\Sp(1)\cdot\Sp(n)$--frame bundle of the quaternionic K\"ahler
 manifold we show that a quaternionic Killing spinor induces a
 Killing spinor on $S$ and a parallel spinor on $\M$. The result
 of M.~Wang then implies that the hyperk\"ahler manifold $\M$ has
 to be locally isometric to $\H^{n+1}$ forcing $M$ to be isometric to the
 quaternionic projective space.

 We would like to thank A.~Swann for several hints and comments
 and W.~Ballmann for encouragement and support.

\section{Semiquaternionic Vector Spaces and Representations}

 The tangent space of a quaternionic
 K\"ahler manifold is not \textit{a priori} a quaternionic left
 vector space because in general the three local complex structures
 are not defined globally. This ambivalence gives rise to the weaker
 notion of a semiquaternionic structure on a real vector space:
 \leer
 \definition
 A semiquaternionic structure on a real vector space $T$ is a
 subalgebra $Q\subset\End\;T$ with $\id_T\in Q$ and $Q\cong\H$
 as $\R$--algebras. It is said to be adapted to a euclidean scalar
 product $\la,\ra$ on $T$ if
 $$
   \la\, q t_1,\,t_2\,\ra = \la\,t_1,\,\overline{q}t_2\,\ra
 $$
 for all $t_1,t_2\in T$ and $q\in Q$, where $\overline{q}$
 denotes conjugation defined by
 $Q\,=\,\R\oplus\mathrm{Im}\,Q\,:=\,\R\;\id_T\oplus [Q,Q]$.
 \leer
 Thus, choosing an isomorphism from $\H$ to $Q$ makes $T$ a quaternionic
 left vector space, however no particular isomorphism is preferred.
 Accordingly, the notion of quaternionic linear map has to be refined:
 \leer
 \definition
 An $\R$--linear map $f:T\to T'$ between vector spaces $T$, $T'$ with
 semiquaternionic structures $Q$, $Q'$ is semilinear, if there exists
 an isomorphism of $\R$--algebras $f^Q: Q\to Q'$ such that $f(qt)=f^Q(q)f(t)$
 for all $t\in T$ and $q\in Q$. If $f$ is semilinear and not identically
 zero $f^Q$ is uniquely defined, because $\id_{T'}\in Q'$ and $Q'\cong\H$
 implies that every non--zero endomorphism in $Q'$ is invertible.
 \leer
 If $T$ is a euclidean vector space with an adapted semiquaternionic
 structure, then the group of all semilinear isometries of $T$ is isomorphic
 to $\Sp(1)\cdot\Sp(n):=\Sp(1)\times\Sp(n)/\Z_2$. Choosing a particular
 isomorphism makes $T$ a true representation of $\Sp(1)\cdot\Sp(n)$ and
 any two representations $T$, $T'$ defined this way are intertwined by a
 semilinear isometry, which is unique up to sign. For this reason we will
 call any such $T$ the defining representation of $\Sp(1)\cdot\Sp(n)$ with
 a choice of isomorphism tacitly understood. It turns out that the defining
 representation $T$ comes along with a preferred isomorphism $\H\to Q$ given
 by the infinitesimal action of $i,\,j,\,k\in\sp(1)\cong\Im\H$ on $T$.

 Similarly, one may construct the defining representation of the
 group $\Sp(n)$ of unitary quaternionic $n\times n$--matrices.
 If $E$ is a complex vector space of dimension $2n$ endowed with a
 symplectic form $\s_E$ and an adapted positive quaternionic structure
 $J$, i.~e.~ $\s_E(\,Je_1,\,Je_2\,)=\overline{\s_E(\,e_1,\,e_2\,)}$
 for all $e_1,e_2\in E$ and $\s_E(\,e,\,Je\,)>0$ for all $0\neq e\in E$,
 then the group of all $\C$--linear symplectic transformations of
 $E$ commuting with $J$ is isomorphic to $\Sp(n)$.
 Choosing a particular isomorphism
 makes $E$ a true $\Sp(n)$--representation and any two representations
 $E$, $E'$ defined this way are intertwined by a $\C$--linear symplectic
 map preserving the quaternionic structure, which is unique up to sign.
 For this reason we will call any such $E$ the defining representation
 of $\Sp(n)$. Note that the Lie algebra of all infinitesimal symplectic
 transformations of $E$ is canonically isomorphic to $\Sym^2E$ with $e_1e_2$
 acting as the endomorphism $\s_E(e_1,\cdot)e_2+\s_E(e_2,\cdot)e_1$, and
 elements of $\Sym^2E$ commute with $J$ if and only if they are real with
 respect to the real structure $\Sym^2J$. Hence, the defining representation
 comes along with a canonical real $\Sp(n)$--equivariant isomorphism
 $\C\otimes_\R\sp(n)\to\Sym^2 E$ of Lie algebras. We will denote the
 defining representation of $\Sp(1)$ by $H$.

 There are several possibilities to give explicit realizations of
 the representations introduced above. In calculations and proofs
 below we will use the following standard picture, differing somewhat
 from Salamon's conventions (cf.~\cite{sala}). Consider the space of
 row vectors $\H^n$ over the quaternions with complex and quaternionic
 structure given by multiplication with $i$ and $j$ from the left. The group
 of unitary quaternionic matrices $\Sp(n):=\{ A\in M_{n,n}\H:\;A^HA=1\}$
 acts on $\H^n$ from the left by multiplying with $A^H$ from the right.
 Thus, it commutes with the complex and quaternionic structure and preserves
 the linear form
 $$
   \s_{\H^n}(\,v_1,\,v_2\,) := [\,v_1v_2^Hj\,]_\C \,,
 $$
 where $[q]_\C:={1\over 2}(q-iqi)\in\C$ is the $\C$--part of $q\in\H$.
 The $\C$--part is obviously $\C$--bilinear, i.~e.~$[xqy]_\C=xy[q]_\C$
 for all $x,y\in\C\subset\H$, and satisfies
 $$
   [\overline{q}]_\C = \overline{[q]}_\C = [-jqj]_\C \,.
 $$
 Using these properties it is easily checked that $\s_{\H^n}$ is indeed
 $\C$--bilinear symplectic and that the quaternionic structure is adapted
 and positive. In this way $\H^n$ becomes the defining representation of
 $\Sp(n)$.

 With a slight modification of the construction above we can make $\H^n$
 the defining representation of $\Sp(1)\cdot\Sp(n)$, too. The
 scalar multiplication with $q\in\H$ from the left on the row
 vectors in $\H^n$ determines a subspace $Q\subset\End\;\H^n$,
 which obviously is a semiquaternionic structure adapted to
 the standard scalar product on $\H^n$ given by
 $$
   \la\,v_1,\,v_2\,\ra:=\mathrm{Re}\; v_1v_2^H =
   \mathrm{Re}\;\s_{\H^n}(\,v_1,\,jv_2\,) \,.
 $$
 The group $\Sp(1)\cdot\Sp(n)=\{ z\cdot A:\;z\in\Sp(1)\textrm{\ and\ }
 A\in\Sp(n)\}$ acts on $\H^n$ from the left through semilinear isometries
 by $(z\cdot A)v:=zvA^H$.

 An important point is particularly obvious in this standard picture and
 in consequence true for every defining representation of $\Sp(1)\cdot\Sp(n)$.
 The infinitesimal action of $i,j,k\in\sp(1)\subset\sp(1)\oplus\sp(n)$ on $T$
 defines a canonical $\Sp(1)\cdot\Sp(n)$--equivariant isomorphism of algebras
 $\H\to Q$ making $T$ a quaternionic left vector space. Thus, the natural
 homomorphism $\Sp(1)\cdot\Sp(n)\to\Aut\;Q$, $z\cdot A\mapsto
 (z\cdot A)^Q$ is trivial on the subgroup $\Sp(n)$ and descends to an
 isomorphism $\SO(3):=\Sp(1)/\Z_2\cong\Aut\;Q$ on the
 complementary subgroup $\Sp(1)$ sending $z\cdot 1\in\Sp(1)$ to
 $z^Q:=(z\cdot 1)^Q$.

 This observation allows a construction which becomes fundamental for
 quaternionic K\"ahler geometry once it is ``gauged''. Let $T$ be the 
 defining representation of $\Sp(1)\cdot\Sp(n)$ and $T'$ an arbitrary
 euclidean vector space with an adapted semiquaternionic structure $Q'$.
 The representations of $\Sp(1)\cdot\Sp(n)$ on $T$ and $\SO(3)$ on $Q$
 define simply transitive right group actions on
 $$
   P=\{ f:\;T\to T'\;\textrm{\ semilinear isometry}\}
 $$
 and
 $$
   S=\{ f^Q:\;Q\to Q'\;\textrm{\ isomorphism of algebras}\}
 $$
 such that the natural projection $f\mapsto f^Q$ is $\Sp(1)$--equivariant.
 Fixing a base point in $P$ to identify $P$ with $\Sp(1)\cdot\Sp(n)$ we get
 the following diagram:
 \begin{center}
 \unitlength1.1pt
 \begin{picture}(160,60)(0,0)
  \put(30,56){$P$}
  \put( 4,2){$S$}
  \put(50,2){$\{T'\}$}
  \put(12,4){\vector(1,0){35}}
  \put(30,53){\vector(-1,-2){21}}
  \put(35,53){\vector(1,-2){21}}
  \put(124,56){$z\cdot A$}
  \put(101,2){$z^Q$}
  \put(155,2){$1$}
  \put(112,4){\vector(1,0){41}}
  \put(130,53){\vector(-1,-2){21}}
  \put(135,53){\vector(1,-2){21}}
 \end{picture}
 \end{center}

 The quaternionic structure on $T$ can be used to construct two important
 relations between the defining representations $H$, $E$ and $T$ of $\Sp(1)$,
 $\Sp(n)$ and $\Sp(1)\cdot\Sp(n)$:

 \begin{Lemma}\label{TisE}
  Let $T$ be the defining representation of $\Sp(1)\cdot\Sp(n)$.
  Define the complex and quaternionic structure on $T$ by the infinitesimal
  action of $i$, $j$ (and $k$) in $\sp(1)\cong\mathrm{Im}\;\H$. With the
  $\C$--bilinear symplectic form $\s_T(\,t_1,\,t_2\,) = \la\,jt_1,\,t_2\,\ra
  + i\la\,kt_1,\,t_2\,\ra$ the vector space $T$ becomes the defining
  representation of $\Sp(n)$. Conversely, the complex and quaternionic
  structure of the defining representation $E$ of $\Sp(n)$ generate a
  subalgebra $Q$ in the $\R$--linear endomorphisms of $E$, which is a
  semiquaternionic structure adapted to the scalar product
  $\la\,\cdot,\,\cdot\,\ra:=\Re\s_E(\,\cdot,\,J\cdot\,)$ making $E$ the
  defining representation of $\Sp(1)\cdot\Sp(n)$. In particular, there
  is up to sign a unique $\Sp(n)$--equivariant, $\C$--linear symplectic
  isomorphism $\Psi:\;T\to E$ preserving the quaternionic structure.
  In the standard picture this isomorphism is simply the identity.
 \end{Lemma}

 This isomorphism has the disadvantage of spoiling the $\Sp(1)$--action on
 $T$. Consequently, it is impossible to use it directly on quaternionic
 K\"ahler manifolds. Nevertheless, we may use it to define a family of
 $\Sp(n)$--equivariant isomorphisms $\C\otimes_\R T\to H\otimes_\C E$
 depending on the choice of a canonical base $p$,\ $q:=Jp$ of $H$ satisfying
 $\s_H(p,q)=1$ by
 \be\label{THE}
  \Phi:\,x\otimes_\R t \longmapsto \tfrac{1}{\sqrt{2}}
  \big(xp\otimes_\C\Psi(t) + xq\otimes_\C J\Psi(t)\big) \,.
 \ee
 All these isomorphisms are real with respect to the real structure
 $J\otimes_\C J$ on $H\otimes_\C E$ and isometries from $\la\,,\,\ra$
 to $\s_H\otimes_\C \s_E$:
 $$
 \begin{array}{rcl}
  \s_H\otimes\s_E(\Phi(1\otimes_\R t_1),\Phi(1\otimes_\R t_2))
  & = & \tfrac{1}{2}
  \big(\s_E(\Psi(t_1),J\Psi(t_2))-\s_E(J\Psi(t_1),\Psi(t_2))\big) \\[1.4ex]
  & = & \tfrac{1}{2}
  \big(\la\,jt_1,\,jt_2\,\ra + i\la\,kt_1,\,jt_2\,\ra
      + \la\,t_1,\,t_2\,\ra - i\la\,it_1,\,t_2\,\ra \big)
  \;\;=\;\; \la\,t_1,\,t_2\,\ra \,,
 \end{array}
 $$
 since $J\Psi(t)=\Psi(jt)$ and $\s_E(\Psi(t_1),\Psi(t_2))=
 \la\,jt_1,\,t_2\,\ra + i\la\,kt_1,\,t_2\,\ra$ by construction. 
 It turns out that there are exactly two canonical bases $p$,\ $q$ for
 which the isomorphism above is not only $\Sp(n)$-- but already
 $\Sp(1)\cdot\Sp(n)$--equivariant, thus defining a more fundamental
 isomorphism better suited for globalization:

 \begin{Lemma}\label{CTisHE}
  Up to sign there is a unique real $\Sp(1)\cdot\Sp(n)$--equivariant
  isomorphism of complex vector spaces
  $$
    \Phi:\quad\C\otimes_\R T \cong H\otimes_\C E \,,
  $$
  which is an isometry from the the $\C$--bilinear extension
  of $\la\,,\,\ra$ to $\s_H\otimes\s_E$.
 \end{Lemma}
 \proof
 It is sufficient to prove this lemma in the standard picture as it
 translates immediately to arbitrary realizations of the defining
 representations. In this picture the canonical base to choose is
 $p=j$ and $q=-1$ (or $p=-j$ and $q=1$) leading to:
 \be\label{CTtoHE}
  \begin{array}{rcccrccc}
   \Phi:\!\! &
    \C\otimes_\R\H^n & \longrightarrow & \H\otimes_\C\H^n &
   \qquad\Phi^{-1}:\!\! &
    \H\otimes_\C\H^n & \longrightarrow & \C\otimes_\R\H^n \\[1.4ex]
   &
    x\otimes_\R v & \longmapsto & {1\over\sqrt{2}}
    (xj\otimes_\C v \, - \, x\otimes_\C jv) &
   &
    q\otimes_\C v & \longmapsto & {1\over\sqrt{2}}
    (1\otimes_\R\overline{q}jv \, + \, i\otimes_\R\overline{q}kv) \,.
  \end{array}
 \ee
 Note that $\Phi^{-1}(iq\otimes_\C v)=\Phi^{-1}(q\otimes_\C iv)=
 i\Phi^{-1}(q\otimes_\C v)$. Thus, $\Phi^{-1}$ is well defined and
 $\C$--linear. As $\Sp(1)\cdot\Sp(n)$--equivariance of $\Phi^{-1}$
 is obvious and $\Phi$ is an isometry it remains to show that $\Phi$
 and $\Phi^{-1}$ are mutual inverses:
 $$
 \begin{array}{rcl}
 (\Phi^{-1}\circ\Phi)(x\otimes_\R v) & = &
 \tfrac{1}{2}\bigl( 1\otimes_\R(-j\overline{x}jv+\overline{x}v)
                + i\otimes_\R(-j\overline{x}ijv+\overline{x}iv)\bigr)\\[1.4ex]
 & = &
 \tfrac{1}{2}\bigl( 1\otimes_\R 2\,\Re\!(x)v
  +i\otimes_\R 2\,\Re\!(-xi)v\bigr)\\[1.4ex]
 & = & x\otimes_\R v \,.
 \end{array}
 $$
 A direct calculation of $\Phi\circ\Phi^{-1}=\id$ is more tedious,
 but can be done. \qed

\section{Principal Bundles on Quaternionic K\"ahler Manifolds}

 ``Gauging'' the pointwise constructions of the previous section leads
 to the definitions of the basic objects of quaternionic K\"ahler geometry.
 However, in strict analogy with K\"ahler geometry one has to impose an
 additional integrability condition:
 \leer
 \definition
 A quaternionic K\"ahler manifold is a Riemannian manifold $M$ of
 dimension $4n$,\ $n\geq 2$ with an adapted semiquaternionic structure
 $Q_xM\subset\End\;T_xM$ on every tangent space which is respected by
 the Levi--Civit\'a connection of $M$:
 $$
   \nabla \G(Q) \subset \G(T^*M\otimes Q)\,.
 $$
 \leer
 Thus, parallel transport of endomorphisms along arbitrary curves $\g$
 induces isomorphisms of $\R$--algebras $Q_{\g(0)}M\to Q_{\g(\t)}M$, and
 \textit{a fortiori} parallel transport of tangent vectors defines
 semilinear isometries $T_{\g(0)}M\to T_{\g(\t)}M$. In particular, the
 Levi--Civit\'a connection is tangent to the reduction of the frame
 bundle to the principal $\Sp(1)\cdot\Sp(n)$--bundle of semilinear
 orthogonal frames
 \be
  P := \{ f: T\to T_xM\;\;\textrm{semilinear isometry}\}
 \ee
 with projection $\pi_M: P\to M,\;f\mapsto x$.
 Additionally, $P$ projects $\Sp(1)$--equivariantly to
 \be\label{s}
  S := \{ f^Q: Q\to Q_xM\;\;\textrm{isomorphism of algebras}\}
 \ee
 via $\pi_S: P\to S,\;f\mapsto f^Q$, which is in turn a principal
 $\SO(3)$--bundle over $M$ with projection
 $\pi: S\to M,\;f^Q\mapsto x$. In this way $P$ may be considered
 as a principal $\Sp(1)\cdot\Sp(n)$--bundle over $M$ or as a principal
 $\Sp(n)$--bundle over $S$:
 \begin{center}
 \unitlength1.2pt
 \begin{picture}(190,60)(0,0)
  \put(30,56){$P$}
  \put( 4,2){$S$}
  \put(52,2){$M$}
  \put(12,4){\vector(1,0){38}}
  \put(30,53){\vector(-1,-2){21}}
  \put(35,53){\vector(1,-2){21}}
  \put(8,30){$\pi_S$}
  \put(49,30){$\pi_M$}
  \put(30,6){$\pi$}
  \put(135,56){$f\!:\!T\!\to\!T_xM$}
  \put(100,2){$f^Q\!:\!Q\!\to\!Q_xM$}
  \put(187,2){$x$}
  \put(153,4){\vector(1,0){32}}
  \put(148,53){\vector(-1,-2){21}}
  \put(165,53){\vector(1,-2){22}}
  \put(126,30){$\pi_S$}
  \put(178,30){$\pi_M$}
  \put(165,6){$\pi$}
 \end{picture}
 \end{center}

 The tangent bundle $TM$ of $M$ is canonically isomorphic to the bundle
 associated to $P$ by the defining representation $T$ sending the
 class $[f,t]\in P\times_{\Sp(1)\cdot\Sp(n)}T$ to $f(t)$. Alternatively,
 this isomorphism can be expressed by the soldering form
 $\theta_M:= f^{-1}\circ(\pi_M)_* \in\G(T^*P\otimes T)$ of $M$.
 Considering the fundamental isomorphism $\Phi:\C\otimes_\R T\to H\otimes E$
 one might try to define vector bundles $\HB$ and $\EB$ from the defining
 representations $H$ and $E$ of $\Sp(1)$ and $\Sp(n)$. In this generality,
 this is only possible on the quaternionic projective spaces, because only
 for these manifolds the bundle $P$ can be covered by a principal
 $\Sp(1)\times\Sp(n)$--bundle. Nevertheless, representations of
 $\Sp(1)\times\Sp(n)$ contained in some $H^{\otimes p}\otimes E^{\otimes q}$
 with $p+q$ even descend to $\Sp(1)\cdot\Sp(n)$ defining vector bundles
 on $M$ associated to $P$. In this way $M$ carries a multitude of naturally
 defined vector bundles besides bundles constructed out of the tangent bundle.
 The usefulness of these vector bundles has been shown by
 \cite{sala} (see also \cite{qklast}).
 In particular, the bundle $\HB\otimes\EB$
 is globally defined and canonically isomorphic to $TM^\C$ as expressed
 e.~g.~by the soldering form
 $\theta_M^{H\otimes E}:=\Phi\circ\theta_M\in\G(T^*P\otimes(H\otimes E))$.

 By definition, the Levi--Civit\'a connection is tangent to $P$ and
 determined by a connection 1--form $\o_M$ with values in the
 direct sum $\,\sp(1)\,\oplus\,\sp(n)$; accordingly, $\o_M$ splits into
 $\,\o_M^{\sp(1)}\oplus\,\o_M^{\sp(n)}$. Additionally, the Levi--Civit\'a
 connection defines a connection on the principal bundle $S$. Recall
 that horizontal lifts $X^h\in T_{f_0}P$ of tangent vectors
 $X={d\over d\t}\big|_0x_\t\in T_{x_0}M$ can be represented by curves
 of semilinear orthogonal frames $f_\t:T\to T_{x_\t}M$ over $x_\t$
 satisfying ${\nabla\over d\t}\big|_0f_\t=0$. Likewise the connection
 on $S$ is defined by representing horizontal lifts $X^h\in T_{f^Q_0}S$
 by curves $f^Q_\t:Q\to Q_{x_\t}M$ of algebra isomorphisms over $x_\t$
 satisfying ${\nabla\over d\t}\big|_0f^Q_\t=0$. More succinctly,
 the connection 1--form $\o$ on $S$ for arbitrary tangent vectors
 ${d\over d\t}\big|_0f^Q_\t$ is given by
 $$
   \o\Big({d\over d\t}\Big|_0f^Q_\t\Big) \; := \;
   \big(f^Q_0\big)^{-1}{\nabla\over d\t}\Big|_0f^Q_\t\,.
 $$
 Using the Leibniz rule for the covariant derivative along curves,
 it is immediately seen that $\pi_S$ projects horizontal vectors
 ${d\over d\t}\big|_0f_\t$ on $P$ to horizontal vectors
 ${d\over d\t}\big|_0f^Q_\t$ on $S$, because for arbitrary $q\in Q,\ t\in T$
 $$
  \Bigl({\nabla\over d\t}\Big|_0f_\t\Bigr)(qt) \; = \;
  \Bigl({\nabla\over d\t}\Big|_0f^Q_\t\Bigr)(q)\,f_0(t)
  \,\,\,+\,\,\,f^Q_0(q)\,\Bigl({\nabla\over d\t}\Big|_0f_\t\Bigr)(t) \,,
 $$
 and consequently ${\nabla\over d\t}\big|_0f_\t=0$ implies
 ${\nabla\over d\t}\big|_0f^Q_\t=0$. Since the projection $\pi_S:P\to S$
 is $\Sp(1)$--equivariant, and vectors tangent to the 
 $\Sp(n)$--action on $P$ are surely vertical with respect to the
 projection $\pi_S$ we conclude:

 \begin{Lemma}\label{piso}
  $\qquad\pi_S^*\o\, = \,\o_M^{\sp(1)}\,.$
 \end{Lemma}

 Remarkably, the curvature of this connection on $S$ depends only on the
 scalar curvature $\k$ of $M$. In fact, according to the classification of
 $\sp(1)\oplus\sp(n)$--curvature tensors due to Alekseevskii
 (cf.~\cite{alex2}, \cite{sala}) the curvature tensor of a quaternionic
 K\"ahler manifold can be expressed in terms of the scalar curvature $\k$
 and a section $\mathfrak{R}$ of $\Sym^4\EB^*$. We have
 \be\label{curv1}
  R = -{\k\over 8n(n+2)}\Bigl(R^H+R^E\Bigr) + R^{hyper},
 \ee
 where $R^H$,$ \; R^E$ and $R^{hyper}$ are $\Sym^2\HB$-- or
 $\Sym^2\EB$--valued 2--forms defined on sections of
 $\HB\otimes\EB\cong TM^\C$:
 \be\label{curv2}
 \begin{array}{rcl}
  R^H_{h_1\otimes e_1,h_2\otimes e_2}
   &=& \s_E(e_1,e_2) h_1\cdot h_2 \in\Sym^2\HB\\[1.4ex]
  R^E_{h_1\otimes e_1,h_2\otimes e_2}
   &=& \s_H(h_1,h_2) e_1\cdot e_2 \in\Sym^2\EB\\[1.4ex]
  R^{hyper}_{h_1\otimes e_1,h_2\otimes e_2} &=&
  \s_H(h_1,h_2) \mathfrak{R}(e_1,e_2,\cdot,\cdot)
  \in\Sym^2\EB^*\cong\Sym^2\EB\,,
 \end{array}
 \ee
 acting as endomorphisms on $\HB\otimes\EB$. Analyzing these terms
 leads to the following formula for the pull--back of the curvature
 2--form $\O$ of the connection $\o$ to $P$:

 \begin{Lemma}\label{curv}
  The curvature 2--form of the connection
  $\o$ on $S$ pulled back to $P$ is given by 
  $$
   \pi_S^*\O\,=\,
   {\k\over 16n(n+2)}\Bigl(
   \la\,\theta_M\wedge i\theta_M\,\ra i
   +\la\,\theta_M\wedge j\theta_M\,\ra j
   +\la\,\theta_M\wedge k\theta_M\,\ra k\Bigr)\,,
  $$
  where by definition $\la\,\theta_M\wedge i\theta_M\,\ra (X,Y):=
  2\la\,\theta_M(X),\,i\theta_M(Y)\,\ra$.
 \end{Lemma}
 \proof
 Instead of calculating the curvature on $S$ directly we will
 calculate the curvature of the bundle $\Sym^2\HB$, which can
 be associated to $S$ or $P$ inheriting the same connection due
 to $\pi_S^*\o=\o_M^{\sp(1)}$. Its curvature considered as an
 $\sp(1)$--valued 2--form on $P$ is thus $\pi_S^*\O$. Obviously
 only $R^H$ acts non--trivially on $\Sym^2\HB$ and neglecting for
 a moment that it is defined for sections of vector bundles we may
 consider it as a real $\Sp(1)\cdot\Sp(n)$--equivariant morphism
 $$
  \C\otimes_\R\L^2T
  \stackrel\cong\longrightarrow\L^2(H\otimes E)
  \stackrel{R^H}\longrightarrow\Sym^2H
  \stackrel\cong\longrightarrow\C\otimes_\R\sp(1) \,,
 $$
 where the first isomorphism is the extension of
 $\C\otimes_\R T\cong H\otimes E$ and the second
 is the canonical isomorphism $\C\otimes_\R\sp(1)\cong\Sym^2H$,
 which makes $H$ the defining representation of $\Sp(1)$. To make
 this isomorphism explicit in the standard picture we choose the
 canonical base $j$,\ $-1$ of $\H$ with $\s_\H(j,-1)=1$ and find
 for the infinitesimal action of $i$,\ $j$ and $k\in\Im\;\H$:
 $$
  \begin{array}{rlcrlcrl}
  i:\!\! & 1 \mapsto -i\,=\, i(-1) &\qquad\qquad &
  j:\!\! & 1 \mapsto -j            &\qquad\qquad &
  k:\!\! & 1 \mapsto -k\,=\, i(-j) \\
   & j \mapsto \phantom{-}k\,=\, i( j) &&
   & j \mapsto \phantom{-}1 &&
   & j \mapsto -i\,=\, i(-1) \,.
  \end{array}
 $$
 Hence, the isomorphism $\C\otimes_\R\sp(1)\cong\Sym^2\H$ maps $i$ to
 $i(1\,j)$,
 $j$ to ${1\over 2}(j^2+1^2)$ and $k$ to ${i\over 2}(j^2-1^2)$. Accordingly,
 the morphism $R^H$ reads in the standard picture 
 $$
 \begin{array}{rcl}
  R^H_{1\otimes_\R v_1,1\otimes_\R v_2}
  & = & \tfrac{1}{2}
   \big(R^H_{j\otimes v_1,j\otimes v_2}-R^H_{1\otimes jv_1,j\otimes v_2}
    -R^H_{j\otimes v_1,1\otimes jv_2}+R^H_{1\otimes jv_1,1\otimes jv_2}
   \big)\\[1.6ex]
  & = & \tfrac{1}{2}
   \big(\s_{\H^n}(v_1,v_2)j^2-\s_{\H^n}(jv_1,v_2)1\,j
    -\s_{\H^n}(v_1,jv_2)1\,j +\s_{\H^n}(jv_1,jv_2)1^2
   \big) \\[1.6ex]
  & = & \mathrm{Re}\;\s_{\H^n}(v_1,v_2)\cdot \tfrac{1}{2}(j^2+1^2)
        +\mathrm{Im}\;\s_{\H^n}(v_1,v_2)\cdot \tfrac{i}{2}(j^2-1^2)
        -\mathrm{Im}\;\s_{\H^n}(v_1,jv_2)\cdot i(1\,j) \\[1.6ex]
  & = & -\la\,v_1,\,jv_2\,\ra j-\la\,v_1,\,kv_2\,\ra k-\la\,v_1,\,iv_2\,\ra i
 \end{array}
 $$
 for $v_1$,\ $v_2\in\H^n$. Being $\Sp(1)\cdot\Sp(n)$--equivariant $R^H$
 can be made a
 $\C\otimes_\R(P\times_{\mathrm{Ad}}\sp(1))\cong\Sym^2\HB$--valued
 2--form in a straightforward way and becomes the $R^H$ defined above.
 Alternatively, $-{\k\over 8n(n+2)}R^H$ can be thought
 of as the horizontal
 $\Sp(1)\cdot\Sp(n)$--equivariant $\sp(1)$--valued curvature form
 $\pi_S^*\O$ on $P$ and is given by the stated formula. The additional
 factor ${1\over 2}$ comes from the definition of the wedge product.
 \qed

\section{The Levi--Civit\'a Connection on $S$ and $\M$}

 Let $\pi: S \rightarrow M $ be the canonical $\SO(3)$--bundle over $M$
 defined in (\ref{s}), with connection form $\o$ and Riemannian
 metric
 $$
  g_S \; = \;  \tfrac{16n(n+2)}{\k} \,
  B(\o, \, \o) \; + \; \pi^* g_M \,,
 $$
 where $B$ is the standard metric on $\sp(1)\cong\Im\H$, and $\k$
 denotes the scalar curvature of $M$. This metric is Einstein, and
 if we rescale the metrics of $M$ and $S$ so that $\k = 16n(n+2)$,
 then $(S,\,g_S)$ has a natural Sasakian 3--structure. The structure
 group of $S$ reduces to $\Sp(n)$, and we can embed the principal
 $\Sp(n)$--bundle $P$ into the frame bundle of $S$ in such a way
 that the soldering form $\theta_S\in\G(T^*P\otimes(\sp(1)\oplus T))$
 of $S$ on $P$ is given by
 \be\label{ts}
  \theta_S \; = \; \sqrt{\tfrac{16n(n+2)}\k} \,
  \o^{\sp(1)}_M \; \oplus \; \theta_M \,.
 \ee
 In this way the Riemannian metric of $S$ is associated to the
 standard metric $B\oplus\la\,,\,\ra$ on $\sp(1)\oplus T$. In terms
 of covariant derivatives the Levi--Civit\'a connection of $g_S$
 is easily computed, and we obtain:
 \begin{Lemma}\label{onie}
  Let $U, \, V$ be vertical vector fields given as fundamental
  vector fields induced by elements of the Lie algebra $\sp(1)$.
  Further, let $X^h, \, Y^h$ be the  horizontal lifts of vector fields
  $X, \, Y$ on $M$ and $\{e_{\nu} \}_{\nu=1,\ldots,4n} $ a locally
  defined horizontal orthonormal frame on $S$. Then the only
  non--zero covariant derivatives are given by
  $$
   \begin{array}{rcl}
   \nabla_U V & = & \tfrac{1}{2} \, [U , \, V] \\[1.4ex]
   \nabla_{X^h} Y^h & = & (\nabla_X Y)^h  \; - \; \tfrac{1}{2} \,
    \O(X^h, \, Y^h) \\[1.4ex]
   \nabla_{X^h} U & = & \nabla_U X^h \;\; = \;\; \tfrac{1}{2}  \, 
    \sum^{4n}_{\nu=1} \, g_S\big( \O(X^h,\,e_{\nu}), \, U\big) \, e_{\nu}\,.
   \end{array}
  $$
 \end{Lemma}
 Here and in the sequel we identify elements of $\sp(1)$ with their
 associated fundamental vector fields.  The Levi--Civit\'a connection
 of $S$ is determined by a $\,\so(\sp(1) \oplus T)$--valued 1--form
 $\o_S$ on the orthonormal frame bundle of $S$, but as the structure
 group reduces to $\Sp(n)$, it is sufficient to know its restriction
 to the $\Sp(n)$--reduction $P$ again denoted by $\o_S$.
 \begin{Lemma}\label{os}
 The connection form $\o_S$ on $P$ can be written as
 $$
 \o_S \; = \; \o^{\sp(n)}_M \; + \; \tfrac{1}{2} \,
 \ad\big(\o^{\sp(1)}_M\big) \; + \; \sqrt{\tfrac\k{16n(n+2)}}  \,
 \big(
 i \,\theta_M \, \wedge \, i \, + \,
 j \,\theta_M \, \wedge \, j \, + \,
 k \,\theta_M \, \wedge \, k 
 \big)\,.
 $$
 \end{Lemma}
 \proof
 For the proof we identify vector fields on $S$ with
 equivariant functions on $P$, i.~e.~
 $$
 \begin{array}{rcl}
 \Gamma(TS)
 & \cong  &
 {\cal C}^\infty(P, \, \sp(1) \oplus T)^{\Sp(n)} \\[1.4ex]
 A & \mapsto  &   \widehat{A}\,.
 \end{array}
 $$
 With respect to this identification the covariant derivative
 translates as
 $$
 \widehat{\nabla_A B}\; = \;d \widehat{B} ( \widetilde{A} ) \; + \;
 \o_S( \widetilde{A}) \, \widehat{B} \,,
 $$
 where $\widetilde{A}$ denotes an arbitrary lift of $A$
 to a vector field on $P$. We will use this formula and
 Lemma \ref{onie} to compute all non--zero terms
 $\o_S( \widetilde{A}) \, \widehat{B}$, which then combine
 to give the stated expression for $\o_S$. The definition
 of the soldering form $\theta_S$ immediately implies that
 the function $\widehat{A}$ is given by $\theta_S(\widetilde{A})$.
 In particular, we have for a fundamental vector field $U$
 \be\label{vert}
 \widehat{U}
 \; = \;
 \theta_S(\widetilde{U} )
 \; = \;
 \sqrt{\tfrac{16n(n+2)}\k} \, \pi^*_S  \o  (\widetilde{U} )
 \; = \;
 \sqrt{\tfrac{16n(n+2)}\k} \, \o  (U )
 \; = \;
 \sqrt{\tfrac{16n(n+2)}\k} \, U \; \in \sp(1) \,.
 \ee

 Let $X^h, \, Y^h$ denote the horizontal lifts of the vector fields
 $ X, \, Y $   on $M$. Because of $ \pi \circ \pi_S = \pi_M $ we
 can assume $ \widetilde{X^h} = \widetilde{X} $. Then,
 $$
 \widehat{ \nabla_{X^h}  Y^h}
  \; =  \; d \widehat{Y}( \widetilde{X})  \; + \;
  \o_S(\widetilde{X}) \,  \widehat{Y} \,.
 $$
 Using Lemma \ref{onie}, Lemma \ref{curv} and equation (\ref{vert}) we find
 $$
 \begin{array}{rcl}
 \o_S(\widetilde{X}) \,  \widehat{Y}
 & = &
 \widehat{(\nabla_X Y)^h}   \; - \;   \tfrac{1}{2} \,
 \widehat{\O( X^h, \, Y^h)} \; - \; d \widehat{Y}( \widetilde{X}) \\[1.6ex]
 & = &
 \o_M (\widetilde{X} ) \,  \widehat{Y}  \; - \;
 \sqrt{\tfrac\k{16n(n+2)}} \,
 \Big(
  \big\la\theta_M(\widetilde{X}),\,i\,\theta_M(\widetilde{Y})\big\ra\, i \;
  +\;\big\la\theta_M(\widetilde{X}),\,j\,\theta_M(\widetilde{Y})\big\ra\, j
  \\[1.6ex]
  &&\qquad\qquad\qquad\qquad
  +\;\big\la\theta_M(\widetilde{X}),\,k\,\theta_M(\widetilde{Y})\big\ra\, k
 \Big)
 \\[1.6ex]
 & = &
 \o_M (\widetilde{X})\,\widehat{Y} \; + \; \sqrt{\tfrac\k{16n(n+2)}}\,
 \Big(
  \big\la i\,\theta_M(\widetilde{X}),\,\widehat{Y}\big\ra \, i \;
  +\;\big\la j\,\theta_M(\widetilde{X}),\,\widehat{Y}\big\ra \, j \;
  +\;\big\la k\,\theta_M(\widetilde{X}),\,\widehat{Y}\big\ra \, k
 \Big)\\[1.6ex]
 & = &
 \o_M (\widetilde{X})\,\widehat{Y} \; + \; \sqrt{\tfrac\k{16n(n+2)}}\,
 \Big(
  i \, \theta_M (\widetilde{X})\, \wedge \, i \; 
  +\;j \, \theta_M (\widetilde{X})\, \wedge \, j \;
  +\;k \, \theta_M (\widetilde{X})\, \wedge \, k
 \Big) \, \widehat{Y} \,.
 \end{array}
 $$
 Let $U, \, V$ be  fundamental vector fields. Then $\widehat{V}$ is
 a constant function and $ d \widehat{V} ( \widetilde{U} ) =
 \widetilde{U}( \widehat{V} ) = 0 $. Hence,
 $$
 \widehat{ \nabla_U V }
 \; = \; \o_S( \widetilde{U} ) \widehat{V}
 \; = \; \tfrac{1}{2} \,  \widehat{[U, \, V]}
 \; = \; \tfrac{1}{2} \, [\widehat{U}, \, \widehat{V}]
 \; = \; \tfrac{1}{2} \, \ad\big(\o_M^{\sp(1)}(\widetilde{U})\big) 
  \,\widehat{V} \,.
 $$
 Finally, let $U$ be fundamental, $X^h$ a horizontal lift of $X$
 and $\{e_{\nu} \}_{\nu = 1, \ldots, 4n} $ a horizontal orthonormal
 frame on $S$, with $E_{\nu} = \theta_S(\widetilde{e_{\nu}}) =
 \widehat{e_{\nu}}$. Then $d\widehat{U}(\widetilde{X}) =
 \widetilde{X} ( \widehat{U} ) = 0 $ and we obtain
 $$
 \begin{array}{rcl}
 \o_S(\widetilde{X}) \, \widehat{U}
 & = &
 \widehat{\nabla_{X^h} U} \; = \; \tfrac{1}{2} \, \sum^{4n}_{\nu =1} \,
 g_S\big( \O(X^h,\,e_{\nu}), \, U \big) \, E_{\nu}\\[1.6ex]
 & = &
 \tfrac\k{16n(n+2)} \, \sum^{4n}_{\nu=1} \,
 \Big(
  \big\la\theta_M(\widetilde{X}),\,i\,E_{\nu}\big\ra\, g_S(i,\,U)\, E_{\nu}\;
  +\;\big\la\theta_M(\widetilde{X}),\,j\,E_{\nu}\big\ra\, g_S(j,\,U)\, E_{\nu}
  \\[1.4ex]
  &&\qquad\qquad\qquad\qquad
  +\;\big\la\theta_M(\widetilde{X}),\,k\,E_{\nu}\big\ra\, g_S(k,\,U)\, E_{\nu}
 \Big)
 \\[1.6ex]
 & = &
 -\,\sum^{4n}_{\nu =1} \,
 \Big(
  \big\la i\,\widehat{X},\,E_{\nu}\big\ra \, B(i,\,U) \, E_{\nu} \;
  +\;\big\la j\,\widehat{X},\,E_{\nu}\big\ra \, B(j,\,U) \, E_{\nu} \;
  +\;\big\la k\,\widehat{X},\,E_{\nu}\big\ra \, B(k,\,U) \, E_{\nu}
 \Big)
 \\[1.6ex]
 & = &
 -\,\sum^{4n}_{\nu =1} \,
 \Big\la\big(i\,B(i,\,U)\;+\;j\,B(j,\,U)\;+\;k\,B(k,\,U)\big)
  \,\widehat{X},\,E_{\nu}
 \Big\ra\,E_{\nu}
 \\[1.6ex]
 & = &
 -\,\o_M^{\sp(1)}(\widetilde{U}) \, \widehat{X} \,.
 \end{array}
 $$
 Combining these three calculations leads to
 $$
 \o_S \; = \; \o_M \;+\; \tfrac{1}{2} \,\ad\big(\o_M^{\sp(1)}\big)
 \;+\; \sqrt{\tfrac\k{16n(n+2)}}\,
 \big(
 i\,\theta_M \,\wedge \, i \,+\,
 j\,\theta_M \,\wedge \, j \,+\,
 k\,\theta_M \,\wedge \, k 
 \big)
 \;-\;
 \o_M^{\sp(1)}\,.
 $$
 In this formula the last summand $-\o_M^{\sp(1)}$ acts by the infinitesimal
 $\sp(1)$--action on $T$ and thus cancels the action of $\o_M^{\sp(1)}$ as
 part of the first summand $\o_M$. So we end up with the stated formula for
 $\o_S$. \qed
 \leer

 Besides the manifold $S$ we also need to consider the cone $\M$
 over $S$, i.~e.~the warped product $\M:=\R^+\times_{t^2} S$
 with metric
 $$
 \widehat{g} \; = \; 
 \tfrac{16n(n+2)}{\k} \, dt^2 \; + \; t^2 \,  g_S .
 $$
 The structure group of $\M$ reduces to $\Sp(n)$, and we can embed
 the principal $\Sp(n)$--bundle $ P_\M:=\R^+\times P$ into the
 bundle of orthonormal frames of $\M$ in such a way that the
 $\R \oplus \sp(1) \oplus T$--valued soldering form on $P_\M$
 is given by
 $$
 \theta_\M
 \; = \;
 \sqrt{\tfrac{16n(n+2)}{\k}} \, (-dt) \; \oplus \; t \, \theta_S
 \; = \;
 \sqrt{\tfrac{16n(n+2)}{\k}} \Big(- dt \; \oplus \; t\, \o_M^{\sp(1)} \Big)
 \; \oplus \; t \, \theta_M  \, .
 $$
 This convention makes the inward pointing vector field $\,\Xi:=
 -\sqrt{\tfrac\k{16n(n+2)}}\,\tfrac{\partial}{\partial t}\,$ correspond
 to $\2:=\theta_\M(\Xi)\in\R$.
 This may be surprising at first but it turns out that only this orientation
 is compatible with a hyperk\"ahler structure on $\M$ introduced
 later. With this choice of soldering form the Riemannian metric $\widehat{g}$
 is associated to the standard metric $\la\,,\,\ra\oplus B\oplus \la\,,\,\ra$
 on $\R\oplus\sp(1)\oplus T$.

 \begin{Lemma}\label{omh}
 The restriction $\o_\M$ of the Levi--Civit\'a connection of
 $\M$ to the reduction $P_\M$ of the bundle of orthonormal
 frames reads:
 $$
 \begin{array}{rcl}
  \o_\M
  & = &
  \o_S \;+\; \sqrt{\tfrac\k{16n(n+2)}} \, \theta_S \, \wedge \, \2 \\[1.6ex]
  & = &
  \o^{\sp(n)}_M \;+\; \tfrac{1}{2} \, \ad\big(\o^{\sp(1)}_M\big)
  \;+\; \o_M^{\sp(1)}\,\wedge\,\2
  \\[1.4ex]
  && +\; \sqrt{\tfrac\k{16n(n+2)}}  \,
  \big(\theta_M \, \wedge \, \2 \,+\,
   i \,\theta_M \, \wedge \, i \, + \,
   j \,\theta_M \, \wedge \, j \, + \,
   k \,\theta_M \, \wedge \, k 
  \big)\,.
 \end{array}
 $$
 In particular, the connection form $\o_\M$ is the pull--back of
 a well defined 1--form on $P$.
 \end{Lemma}
 \proof
 The proof is similar to the proof of the corresponding formula
 for $\o_S $. Let $ \widehat{\nabla}$ denote the covariant
 derivative for the  Levi--Civit\'a connection  of $\widehat{g}$.
 The only non--vanishing terms are
 $$
 \widehat{\nabla}_X Y
 \;=\; \nabla^S_X Y \;+\; \sqrt{\tfrac\k{16n(n+2)}} \, \tfrac{1}{t} \, 
 \widehat{g} (X, \, Y) \, \Xi
 \qquad\mbox{and}\qquad
 \widehat{\nabla}_X \Xi
 \;=\; \widehat{\nabla}_\Xi X
 \;=\; -\sqrt{\tfrac\k{16n(n+2)}} \, \tfrac{1}{t} \, X \,.
 $$
 Using the same notation as in the proof of Lemma \ref{os} we obtain
 $$
 \o_\M(\widetilde{X})\widehat{Y}
 \;=\;
 \o_{S}(\widetilde{X})\widehat{Y}
 \;+\; \sqrt{\tfrac\k{16n(n+2)}} \, \tfrac{1}{t} \,
 \big\la
  \theta_\M(\widetilde{X}), \,
  \theta_\M(\widetilde{Y})
 \big\ra \,\widehat{\Xi}
 \;=\;
 \o_{S}(\widetilde{X})\widehat{Y}
 \;+\; \sqrt{\tfrac\k{16n(n+2)}} \,
 \big\la\theta_{S}(\widetilde{X}),\,\widehat{Y}\big\ra \,\2 \,.
 $$
 Having this expression for horizontal vector fields
 $X,\,Y$ we immediately derive
 $$
 \o_\M 
 \;=\; \o_S \;+\; \sqrt{\tfrac\k{16n(n+2)}}\,\theta_S \, \wedge \, \2
 \;=\; \o_S \;+\; \sqrt{\tfrac\k{16n(n+2)}}\,\theta_M \, \wedge \, \2
 \;+\; \o_M^{\sp(1)} \, \wedge \, \2 \,.\qed
 $$
 \leer
 An interesting application of the formulas given above relates the curvature
 tensor $\widehat{R}$ of the manifold $(\M,\,\widehat{g})$ to the
 hyperk\"ahler part $R^{hyper}$ of the curvature tensor of $(M,\,g)$. The
 proof will be given in appendix \ref{proofrhyper}.
 \begin{Proposition}\label{hypcurv}
 The curvature tensor $\widehat{R}$ is horizontal with respect to
 $\widehat{\pi}: \M\to M$. Its only non--vanishing terms are
 $$
 \widehat{R}(X,\,Y) \;=\; R^{hyper}(\widehat{\pi}_*X,\,\widehat{\pi}_*Y) \,.
 $$
 The right--hand side acts on the orthogonal complement
 $(T^V\M)^\perp$ of the vertical tangent bundle,
 which is canonical isomorphic to $\widehat\pi^* TM$.
 \end{Proposition}
 The manifold $\M$ has been previously studied by A.~Swann
 using the notation $\mathcal{U}(M)$ \cite{swann}. He constructs $\M$
 as the $\Z_2$--quotient of the total space of the locally defined bundle
 $\HB$ with zero section removed. The metric $\widehat{g}$ is a member of
 the family of hyperk\"ahler metrics on $\M$ introduced in
 \cite{swann}. In particular, the proposition above is implicit in his work
 (see also \cite{swann1}).

\section{Quaternionic Killing Spinors}

 In this chapter we recall the quaternionic Killing equation introduced
 in \cite{qklim}. It will provide us with an equivalent formulation of the 
 limit case. First, we have to collect some facts for the spinor bundle
 and  Clifford multiplication on a quaternionic K\"ahler manifold
 (cf.~\cite{qklast}). 

 The spinor bundle of a $4n$--dimensional quaternionic K\"ahler spin
 manifold decomposes into a sum of $n+1$ subbundles, which can be
 expressed using the locally defined bundles $\EB$ and $\HB$
 (cf.~\cite{barsal}, \cite{hijmil} or \cite{wang}). For this we have
 to introduce the bundles $ \L^s_\circ \EB $. They are associated to
 the irreducible $\Sp(n)$--representations on the spaces $\L^s_\circ E$
 of primitive $s$--vectors, which are the kernels of the contraction
 $\s_E\es : \L^sE \longrightarrow \L^{s-2}E$ with the symplectic form
 $\s_E$. With this notation the spinor bundle can be written
 \be\label{spinor}
 \spb(M)
 \; =\; \bigoplus^n_{r=0}\,\spb_r(M)
 \;:=\; \bigoplus^n_{r=0}\,\Sym^r\HB \otimes \L^{n-r}_\circ\EB \, .
 \ee
 In order to define the Clifford multiplication we have to fix notations
 for modified contraction and multiplication on $\Sym^r \HB$ and
 $\L^{s}_\circ\EB$. Contraction preserves the primitive spaces,
 i.~e.~if $\eta$ is in $\L^s_\circ E $ then $e^\#\,\es\,\eta\;\in\,
 \L^{s-1}_\circ E$, where $e^\#:=\s_E(e,\,\cdot)\in E^*$ denotes
 the dual of $e\in E$. However, this is not true for the wedge product
 and the projection $e\dn\eta$ of $e\wedge\eta$ onto $\;\L^{s+1}_\circ E$
 is given by
 $$
  e \dn \eta =  e \wedge \eta  \;  - \;
  \tfrac{1}{n-s+1} \, L_E \wedge (e^\#\,\es\,\eta) \, ,
 $$
 where $L_E$ is the canonical bivector associated to $\s_E$
 under the isomorphism $\L^2 E \cong \L^2 E^*$. Let $h\cdot$
 denote the symmetric product with $h\in H$, and for
 $h^\#:=\s_H(h,\cdot)\in H^*$ we define $h^\#\,\esn: \Sym^r H
 \rightarrow \Sym^{r-1} H$ by $h^\#\,\esn:=\tfrac{1}{r} h^\#\,\es$.
 Let $h\otimes e\in\HB\otimes\EB = TM^{\C}$ be a tangent vector.
 Then, the Clifford multiplication
 $\mu( h \otimes e) : \spb(M) \rightarrow \spb(M)  $ is given by
 \be\label{mult}
  \mu(h \otimes e) \quad = \quad \sqrt{2} \;
  (h \cdot \,  \otimes\, e^\# \, \es \,
  \; + \;
  h^\# \, \esn \, \otimes e \, \dn  \,)\, .
 \ee
 In particular, it maps the subbundle
 $\spb_r(M)$ to the sum $\spb_{r-1}(M) \oplus \spb_{r+1}(M)$ and thus
 splits into two components
 $$
 \mu^+_- : \quad  TM \,\otimes\,\spb_r(M)\; \longrightarrow\;\spb_{r+1}(M)
 \quad\mbox{and}\quad
 \mu^-_+ : \quad  TM\,\otimes\,\spb_r(M)\; 
 \longrightarrow \; \spb_{r-1}(M)  \, ,
 $$
 with $\;\mu^+_-(e\,\otimes\,h)=\sqrt{2} \,(h\cdot\otimes\; e^\#\,\es \;)$
 and  $\;\mu^-_+(e\,\otimes\,h)=\sqrt{2} \,(h^\#\,\esn\,\otimes e\,\dn \;)$.
 There are two operations defined similar  to Clifford multiplication
 $$
 \begin{array}{rccc}
  \mu^+_+: & TM \otimes \Sym^r \HB \otimes \L^s_\circ\EB & \longrightarrow &
   \Sym^{r+1} \HB \otimes \L^{s+1}_\circ\EB \\[1.4ex]
  & h \otimes e \otimes \psi & \longmapsto &
   {\sqrt 2}\, (h \cdot \, \otimes \,  e \, \dn\,)\psi
 \end{array}
 $$
 and
 $$
 \begin{array}{rccc}
  \mu^-_-: & TM \otimes \Sym^r \HB \otimes \L^s_\circ\EB & \longrightarrow &
   \Sym^{r-1} \HB \otimes \L^{s-1}_\circ\EB \\[1.4ex]
  & h \otimes e \otimes \psi &\longmapsto &
   {\sqrt 2}\, (h^\#\,\esn\,\otimes\, e^\#\,\es\,)\psi \, .
 \end{array}
 $$
 Using these notations the Dirac operator $D$ can be written $D=D^+_- +D^-_+$
 with
 $$
  D^{+}_{-} \,:=\,
  \mu^{+}_{-} \circ \nabla : \spb_r(M)\,\longrightarrow\,\spb_{r +1}(M)\,
  \qquad 
  D^{-}_{+} \,:=\,
  \mu^{-}_{+} \circ \nabla : \spb_r(M)\,\longrightarrow\,\spb_{r-1}(M)\, .
 $$
 The square of the Dirac operator respects the splitting of the spinor bundle,
 i.~e.~$D^2:\spb_r(M)\,\longrightarrow\,\spb_r(M)$ and we have $ D^+_- D^+_- =
 0 = D^-_+  D^-_+ $.
 \leer

 A quaternionic Killing spinor is by definition (cf.~\cite{qklim}) a
 section $\psi = (\psi_0, \, \psi_1,\, \psi_-) $ of the Killing bundle
 $$ 
 \spb^{Killing}(M)
  \;:=\; \spb_0(M)\,\oplus\,\spb_1(M)\,\oplus\,\L^{n-2}_\circ\EB
  \;\cong\;
  \L^n_\circ\EB\,\oplus\,(\HB\otimes\L^{n-1}_\circ\EB)
  \,\oplus\,\L^{n-2}_\circ\EB\, ,
 $$
 satisfying the following \textit{quaternionic Killing equation}
 for some parameter $\l\ne 0$ and all tangent vectors $X$:
 $$
 \begin{array}{rcl}
  \nabla_X\psi_{0} & = & \phantom{-\,\tfrac\l{4n} \,\mu^+_-(X)\,\psi_0}
  \,-\,\tfrac\l{n+3}\,\mu^-_+(X)\,\psi_1\\[1.4ex]
  \nabla_X\psi_{1} & = &  -\,\tfrac\l{4n} \,\mu^+_-(X)\,\psi_0
  \phantom{\,-\,\tfrac\l{n+3}\,\mu^-_+(X)\,\psi_1}
           \,+\,\tfrac{3\l}{2(n+3)}\,\mu^+_+(X)\,\psi_-\\[1.4ex]
  \nabla_X\psi_{-} & = & \phantom{-\,\tfrac\l{4n} \,\mu^+_-(X)\,\psi_0}
  \,-\,\tfrac\l{4n} \,\mu^-_-(X)\,\psi_1 \,.
 \end{array}
 $$
 We remark that if $ (\psi_0, \, \psi_1,\, \psi_-) $ is a solution
 for parameter $\l$, then $(\psi_0, \, -\psi_1,\, \psi_-)$
 is a solution for parameter $-\l$. In
 \cite{qklim} we showed that for any solution $\psi\ne 0$ the spinor
 $\psi_0+\psi_1$ is an eigenspinor for the minimal eigenvalue
 $$
  \l\; = \;\pm\sqrt{ \tfrac\k{4} \, \tfrac{n+3}{n+2}} \,.
 $$
 In particular, only for two values of the parameter $\l\neq 0$
 there can possibly exist non--trivial solutions. Conversely, any
 eigenspinor for the minimal eigenvalue is of the form $\psi_0+\psi_1
 \in\G(\spb_0(M)\oplus \spb_1(M))$, and the augmented eigenspinor
 $(\psi_0,\psi_1,\psi_-)$ with
 $\psi_-:=\tfrac{1}{4\l}\,\tfrac{n+3}{n+4} \,(\mu^-_-\circ \nabla)\psi_1
 \in\G(\L^{n-2}_\circ\EB)$ is a solution of the quaternionic Killing
 equation. Obviously, solutions are sections parallel with respect to
 a modified connection. Its curvature is precisely the hyperk\"ahler part
 $R^{hyper}$ of the curvature of $M$. Since this part vanishes on the
 quaternionic projective space the Killing bundle of $\HP n$ is trivialized
 by augmented eigenspinors with minimal eigenvalue. Proposition \ref{hypcurv}
 then motivates to lift a quaternionic Killing spinor to a parallel spinor
 on $\M$.

 For our purpose of characterizing the limit case it is more convenient to
 consider an equivalent version of the quaternionic Killing equation. Let
 $\psi = (\psi_0, \, \psi_1,\, \psi_-) $ be a solution of the original 
 equation, then the scaled section 
 $$
 \psi^{scal} 
 \; := \;
 (\psi^{scal}_0, \, \psi^{scal}_1,\, \psi^{scal}_-)
 \; = \;
 \Big(\sqrt{\tfrac{n+3}{4n}} \, \psi_0, \;  \psi_1,\,
 - \, \sqrt{\tfrac{4n}{n+3}} \, \psi_-\Big) 
 $$
 is a solution of the equation
 $$
 \nabla_X \, \psi^{scal} \; = \; 
  -\sqrt{\tfrac{\k }{16n(n+2)}} \, A_X \,\psi^{scal} \,,
 $$
 where $ \psi^{scal} $ is considered as column vector
 with three entries and $A_X$ denotes the matrix
 $$
 A_X \;=\;\left(\begin{array}{ccc}
 0 &  \mu^-_+(X) & 0 \\[1.6ex]
 \mu^+_-(X) & 0 & \tfrac{3}{2}\,\mu^+_+(X) \\[1.6ex]
 0 & - \, \mu^-_-(X) & 0 
 \end{array}\right)\,.
 $$
 For the remainder of this article, the index denoting the scaling
 will be omitted. The following lemma shows that $A_X$
 can be interpreted as part of an $\sp(n+1)$--action.
 \begin{Lemma}\label{decomp}
 Let $ F = H \oplus E$ be the defining representation of $\Sp(n+1)$
 with symplectic form $\s_F=\s_H+\s_E$. Restricted to the subgroup
 $\Sp(1)\times\Sp(n)$ the $\Sp(n+1)$--representation $\L^s_\circ F$
 decomposes into 
 $$
 \L^s_\circ F \cong
 \L^s_\circ E \, \oplus \,(H \otimes \L^{s-1}_\circ E)
 \,\oplus \, \L^{s-2}_\circ E \,,
 $$
 which descends to a well defined representation of $\Sp(1)\cdot\Sp(n)$
 if $s$ is even. Explicitly, the isomorphism 
 $
 \iota : \L^s_\circ E \, \oplus \, (H \otimes \L^{s-1}_\circ E)
 \,\oplus \, \L^{s-2}_\circ E \longrightarrow \L^s_\circ F $
 is given by 
 $$
 \iota \big( \phi_0 \oplus (h \otimes \phi_1) \oplus \phi_- \big)
 \; =  \;
 \phi_0 \,+\, (h \wedge \phi_1) \,+\,
 (L_H - \tfrac{1}{n-s+2}L_E) \wedge \phi_-  \,.
 $$
 Similarly, $\,\Sym^2 F \cong \C\otimes_\R\sp(n+1)$ decomposes into
 $ \Sym^2 H \, \oplus \, (H \otimes E) \, \oplus \, \Sym^2 E$.
 For $s=n$, the subspace $ H \otimes E \subset \Sym^2 F $ acts on
 $\L^n_\circ F$ via
 $(h \otimes e)\,\phi \;=\; \frac{1}{\sqrt{2}} \, A_{h \otimes e} \, \phi$.
 \end{Lemma}
 \proof
 It is clear that $\iota$ is an injective map to $\L^s F$. It remains
 to show that its image is already contained in $\L^s_\circ F$. Since
 $\s_F = \s_H + \s_E$ the statement follows from
 $$
 \s_F \, \es \, (L_H - \tfrac{1}{n-s+2}\, L_E) \, \wedge \, \phi_-
 \; = \;
 [\s_H \, \es, \, L_H \wedge] \, \phi_- \; - \; \tfrac{1}{n-s+2} \,
 [\s_E \, \es, \, L_E \wedge] \, \phi_-
 \; = \; 0 \,,
 $$
 where we used the relations $[\s_H\,\es,\,L_H\wedge]\,\phi_-\,=\,\phi_-$
 and $[\s_E\,\es,\,L_E\wedge]\,\phi_- \, = \,(n-s+2)\,\phi_- $. Comparing
 dimensions shows that $ \iota $ is in addition surjective, hence it defines
 an isomorphism.

 The action of an element $ f_1 \, f_2 \in \Sym^2 F $ on $F$ is given by 
 $(f_1 \, f_2)(f) = \s_F(f_1 ,\, f) f_2 \, + \,  \s_F(f_2 ,\, f) \, f_1 \,$.
 It extends as derivation to $ \L^*_\circ F $ and can be explicitly written as
 $(f_1 \, f_2)(\o)\,=\, (f_2 \wedge f_1^\# \es \, + \, f_1 \wedge f_2^\# \es )
 (\o)$. Hence, for $h \in H$ and $ e \in E$ considered as elements of $F$,
 the element $h\otimes e\in H\otimes E$ is identified with $h\,e\in\Sym^2 F$,
 and the action on $\iota(\phi) = \phi_0 \oplus (a \wedge \phi_1) \oplus
 (L_H - \tfrac{1}{2}L_E) \wedge \phi_-\in \L^n_\circ F$ is given by
 $$
 \begin{array}{rcl}
 (h \otimes e) \, \iota(\phi_0)
 &=&
 h \wedge e^\# \es \phi_0 \\[1.4ex]
 &=&
 \frac{1}{\sqrt{2}}\,\iota\big(\,\mu^+_-(h \otimes e) \phi_0\,\big)
 \\\\
 (h \otimes e)\, \iota(a\otimes\phi_1)
 &=&
 (h \wedge e^\# \es \; + \; e \wedge h^\# \es)(a \wedge \phi_1) \\[1.4ex]
 &=&
 - \; h \wedge a \wedge e^\# \es  \, \phi_1
 \; + \;
 \s_H(h, \, a) \, e \wedge \, \phi_1 \\[1.4ex]
 &=&
 -\s_H(h, \, a) \, \big(L_H \; - \; \tfrac{1}{2} L_E\big)
 \, \wedge \, e^\# \es \, \phi_1
 \; + \; \s_H(h, \, a) \, e \dn \, \phi_1 \\[1.4ex]
 &=&
 \frac{1}{\sqrt{2}}\,\iota
 \Big(
 \mu^-_+ (h \otimes e)  \, (a \otimes \phi_1)
 \; - \;
 \mu^-_- (h \otimes e) \, (a \otimes \phi_1)
 \Big)
 \end{array}
 $$
 and
 $$
 \begin{array}{rcl}
 (h \otimes e)\, \iota(\phi_-)
 &=&
 - \, \tfrac{1}{2} \, h \wedge e^\# \es (L_E \wedge \phi_-)
 \; + \;
 e \wedge h^\# \es (L_H \wedge \phi_-) \\[1.4ex]
 &=&
 \tfrac{1}{2} \, h  \wedge e \wedge \phi_-
 \; - \;
 \tfrac{1}{2} \, h \wedge L_E \wedge (e^\# \es \, \phi_-)
 \; - \;
 e \wedge h \wedge \phi_- \\[1.4ex]
 &=&
 \tfrac{3}{2}\, h \wedge e\dn \phi_- \\[1.4ex]
 &=&
 \tfrac{3}{2\sqrt{2}} \, \iota\big(\,
     \mu^+_+(h \otimes e) \phi_- \,\big) \,.
 \end{array}
 $$
 Hence, we see that operation of $(h\otimes e)$ on $\iota(\phi)$ is just
 application of the matrix $\frac{1}{\sqrt{2}}\,A_{h\otimes e}$ to the
 column vector $\phi$.
 \qed
 \leer
 To stress the origin of the operation $\frac{1}{\sqrt{2}}A_X$ from a
 group action,
 we introduce $\star : H \otimes E \otimes \L^n_\circ F \rightarrow
 \L^n_\circ F $ for the infinitesimal action of $H\otimes E$ on 
 $\L^n_\circ F$. With this notation the quaternionic Killing equation
 reads
 \be\label{qkgl}
  \nabla_X \psi 
  \;=\; - \, \sqrt{\tfrac\k{16n(n+2)}} \, A_X \psi 
  \;=\; - \, \sqrt{\tfrac\k{ 8n(n+2)}} \, \Phi(X) \star \psi \,,
 \ee
 where $\Phi$ is the isomorphism defined in Lemma \ref{CTisHE}.

\section{The Geometry of $\M$ and Application to Spinors}

 The aim of this section is to show that the quaternionic Killing equation,
 considered as a differential equation on equivariant functions on $P$ can
 be interpreted in three different ways: first, of course, when we think of
 its solutions as sections of the Killing--bundle $\spb^{Killing}(M)$ on $M$
 they are quaternionic Killing spinors. Interpreted as a section of the spinor
 bundle on $S$, the solutions are Killing spinors, and finally solutions
 pulled back to $P_\M=\R^+\times P$ are parallel sections of the
 spinor bundle of $\M$.

\subsection{The Hyperk\"ahler Structure of $\M$}
\label{sixone}

 First we recall that the structure group of both $S$ and $\M$
 reduce to $\Sp(n)$, i.~e.~the tangent bundles may be associated to the
 principal $\Sp(n)$--bundles $P$ and $P_\M$ through the
 $\Sp(n)$--representations $\sp(1)\oplus T$ and $\R\oplus\sp(1)\oplus T$
 respectively. However, as $\Sp(n)$--representation $T$ can be identified
 with $E$ according to Lemma \ref{TisE} reflecting the isomorphism
 $$
  \pi^* TM \cong P\times_{\Sp(n)}T \cong P\times_{\Sp(n)}E
 $$
 of vector bundles on $S$. More important, this identification is
 respected by the Levi--Civit\'a connection, because the infinitesimal
 $\sp(1)$--action on $T$ present in $\o_M=\o_M^{\sp(1)}\oplus\o_M^{\sp(n)}$
 is canceled in the connection form $\o_S$ of $S$.
 For this reason we will consequently identify $T$ with $E$ and consider
 $E$ as a euclidean vector space with scalar product
 $\la\,\cdot,\,\cdot\,\ra=\Re\s_E(\cdot,J\cdot)$.

 In the same vein we combine the obvious identification
 $\R\oplus\sp(1)\cong\R\oplus\Im\;\H=\H$ with the isomorphism $\H\to H$
 of defining representations of $\Sp(1)$ unique up to sign to get an
 isometry from the standard metric on $\R\oplus\sp(1)$ to $H$ with
 scalar product $\Re\s_H(\cdot,J\cdot)$ sending $\2$, $i$, $j$, $k$
 to $\1$, $\I$, $\J$, $\K\in H$. In this way we get an isometry
 $$
  \R\oplus\sp(1)\oplus T \longrightarrow F
 $$
 with the defining representation $F=H\oplus E$ of $\Sp(n+1)$ considered
 as a euclidean vector space with scalar product $\Re\s_F(\cdot,J\cdot)$.
 As this isometry is $\Sp(n)$--equivariant by construction the tangent
 bundles of $\M$ and $S$ are associated to the
 $\Sp(n)$--representations $F$ and $\{\1\}^\perp=:(\Im\,H)\oplus E\subset F$.

 Though the structure group of $\M$ reduces to $\Sp(n)$, the
 holonomy of $\M$ does not. Nevertheless, we will show that it
 is a subgroup of $\Sp(n+1)$, i.~e.~$\M$ is hyperk\"ahler. For
 this purpose we group the summands of the connection form $\o_\M$
 of $\M$ given in Lemma \ref{omh} as follows
 $$
 \begin{array}{rcl}
  \o_\M^{\sp(1)} &:=& \frac{1}{2}\,\ad(\o_M^{\sp(1)})
    \;+\; \o_M^{\sp(1)} \wedge \1 \\[1.4ex]
  \o_\M^{\sp(n)} &:=& \o_{M}^{\sp(n)} \\[1.4ex]
  \o_\M^{\H^n} &:=& \sqrt{\frac{\k}{16n(n+2)}}\,
   \big(\theta_M\wedge\1 + i\theta_M \wedge \I
   + j\theta_M \wedge \J + k\theta_M \wedge \K \big)\,.
 \end{array}
 $$
 Thus, the Levi--Civit\'a connection can be written
 $\o_\M= \o_\M^{\sp(1)} \,+\,
 \o_\M^{\sp(n)} \,+\, \o_\M^{\H^n}$.

 \begin{Lemma}\label{osp1}
 The actions of $\o_\M^{\sp(1)}$ and $\o_M^{\sp(1)}$
 are the same, i.~e.~
 \be
  \o_\M=\o_M^{\sp(1)} \,+\, \o_M^{\sp(n)}
  \,+\,\sqrt{\tfrac\k{16n(n+2)}}
   \,\big(
    \theta_M\wedge\1 + i\theta_M\wedge\I
    + j\theta_M\wedge\J + k\theta_M\wedge\K
   \big)\,.
 \ee 
 In particular, the connection form $ \o_\M$ takes values in
 $\sp(n+1)$.
 \end{Lemma}
 \proof
 Note that for $q\in\H$ and imaginary $z\in\Im\;\H$ we have
 $\frac{1}{2}(zq+qz)=(\Re\;q)z-\la\,q,\,z\,\ra=-(z\wedge 1)q$.
 Using this algebraic identity the infinitesimal $\sp(1)$--action
 on $\H$ in the standard picture can be written as
 $$ 
  -qz = \tfrac{1}{2}\ad(z)q + (z\wedge 1)q \,.
 $$
 Hence, a particular merit of the identifications above is that the two
 summands $\tfrac{1}{2}\ad(\o_M^{\sp(1)})$ and $\o_M^{\sp(1)}\wedge\1$ of
 the Levi--Civit\'a connection of $\M$ acting on $\R\oplus\sp(1)$
 combine into the infinitesimal action of $\o_M^{sp(1)}$ on $H$.
 Consequently, the summands $\o_\M^{\sp(1)}$ and
 $\o_\M^{\sp(n)}$ take values in $\sp(n+1)$, i.~e.~in
 the infinitesimal quaternionic linear isometries of $F$. The same
 is true for $\o_\M^{\H^n}$ because of its $\H$--linearity.
 \qed
 \leer

 With the Levi--Civit\'a connection being $\sp(n+1)$--valued the manifold
 $\M$ is hyperk\"ahler, and we may use the description of the spinor
 bundle for the more general quaternionic K\"ahler manifolds given in
 (\ref{spinor}). Consider a complex vector space $\C^2$ endowed
 with a symplectic form $\s_{\C^2}$ and a positive quaternionic
 structure $J$.
 Choosing an isomorphism to the group of all symplectic transformations
 of $\C^2$ commuting with $J$ would make $\C^2$ the defining representation
 of $\Sp(1)$.
 However, on a hyperk\"ahler manifold this $\Sp(1)$--symmetry is not a local
 ``gauged'' symmetry, but a purely global one. For this reason the trivial
 $\C^2$--bundle on $\M$ plays the role of $H$ on $M$. Accordingly,
 the spinor bundle of $\M$ is associated to the
 $\Sp(n)$--representation
 $$
  \S \,=\,\bigoplus_{r=0}^{n+1}\S_r
  \,=\,\bigoplus_{r=0}^{n+1}\Sym^r\C^2\otimes\L_\circ^{n+1-r} F\,,
 $$
 where $\Sp(n)$ operates trivially on $H \subset F$.
 The Clifford multiplication with complex tangent vectors in
 $\C^2\otimes F$ is then given by the formula (\ref{mult}). To
 describe the Clifford multiplication with real tangent vectors
 however, we have to choose an isomorphism among the family of isometries
 $\Phi:\C\otimes_\R F\to\C^2\otimes F$ defined in equation (\ref{THE}).
 For quaternionic K\"ahler manifolds this isometry is uniquely fixed by
 the additional local $\Sp(1)$--symmetry up to sign, but this is no longer
 true in the hyperk\"ahler case. In fact, we get a family of Clifford
 multiplications depending on the choice of $\Phi$, i.~e.~of a canonical
 base $p$,\ $q$ of $\C^2$ satisfying $Jp=q$ and $\s_{\C^2}(\,p,\,q\,)=1$.
 All these are intertwined by the global $\Sp(1)$--symmetry acting on
 $\C^2$. In this sense, to define the Clifford multiplication with real
 tangent vectors $f\in F$, we first have to apply the isomorphism
 \be\label{FtoF}
  \Phi: 1\otimes_{\R}f \longmapsto
   \tfrac{1}{\sqrt{2}}\big( p\otimes f - q\otimes Jf \big) \,.
 \ee
 Note that this isomorphism points out the ambivalence of the vector bundle
 associated to the defining representation $F$ of $\Sp(n+1)$. Normally, this
 vector bundle is the real tangent bundle, but in the description of the
 spinor bundle it plays a role strictly analogous to the isotropic subspace
 $T^{0,1}M$ of the complexified tangent bundle on K\"ahler manifolds.

 To describe the spinor bundle on $S$ we recall that the spinor module
 of $\Spin(4n+4)$ decomposes into the two half--spin modules $\S^\pm$.
 When restricted to $\Sp(n+1)$ the representations $\S^\pm$ decompose
 further into a sum of certain $\S_r$. Since Clifford multiplication
 maps $\S_r$ to $\S_{r-1}\oplus\S_{r+1}$ and interchanges
 $\S^+$ with $\S^-$, we conclude
 $$
 \S^+\,=\,\bigoplus_{{r=0\atop r\equiv 1\,(2)}}^{n+1}
  \Sym^r\C^2\otimes\L_\circ^{n+1-r} F
 \qquad\qquad
 \S^-\,=\,\bigoplus_{{r=0\atop r\equiv 0\,(2)}}^{n+1}
  \Sym^r\C^2\otimes\L_\circ^{n+1-r} F\,,
 $$
 at least if $n$ is even. Due to Lemma \ref{decomp} $\S^+$ and $\S^-$
 are equivalent as $\Sp(n)$--representations, and the spinor module of $S$
 is associated to the $\Sp(n)$--principal bundle $P$ through either one,
 e.~g.~$\S^+$. The Clifford multiplication with $f\in\{\1\}^\perp
 \subset F$ is then given by 
 $f\mathop{\cdot}\limits_S:=(f \wedge \1)\cdot= f \cdot \1 \cdot$.

\subsection{Reinterpretation of the Quaternionic Killing Equation}

 In this section we will translate the quaternionic Killing equation
 (\ref{qkgl}) on $M$ into the equation for a parallel spinor on $\M$. 
 Let $\psi$ be a quaternionic Killing spinor on $M$, which we will consider
 as an equivariant function on $P$, i.~e.~$\psi\in{\cal C}^\infty(P,\,
 \L_\circ^nF)^{\Sp(n)}$. Then the quaternionic Killing equation reads
 \be \label{qkfgl}
  \big(d \,+\, \o_M 
   \,+\, \sqrt{\tfrac\k{8n(n+2)}}\,\theta_M^{H\otimes E}\star
  \big)\psi \;=\; 0\,,
 \ee
 where by definition $\theta_M^{H\otimes E}:= \Phi \circ \theta_M$.
 The subbundle $\spb_1(\M)$ of the spinor bundle of
 $\M$ is associated to $P_\M$ by the representation
 $\C^2\otimes\L_\circ^nF$. To construct a section of this bundle, we
 proceed as follows: we lift $\psi$ to $P_\M$ by extending it
 constantly along the $\R^+$--direction and choose an arbitrary constant
 vector $\xi\in \C^2$. Then, the function $\xi\otimes \psi$ defines a
 section in $\spb_1(\M)$. As function on $P_\M$ it
 satisfies
 $$
 \big(
  d \,+\, \o_M \,+\, \sqrt{\tfrac\k{8n(n+2)}}
  \,\id\otimes\theta_M^{H\otimes E}\star
 \big)(\xi\otimes\psi) \;=\; 0\,,
 $$
 where the connection and soldering form are considered to be pulled back to 
 $P_\M$. Obviously, they do not act on $\xi$ but only on $\psi$.
 Using the crucial Lemma \ref{thetastar} below we can replace
 $\sqrt{\tfrac\k{8n(n+2)}}\,\id\otimes\theta_M^{H\otimes E}\star\,$ by
 $\,\o_\M^{\H^n}$, so that the resulting equation on 
 $P_\M$ reads
 $$
 \big(d \,+\, \o_M^{\sp(1)}
   \,+\, \o_M^{\sp(n)} \,+\, \o_\M^{\H^n}
 \big)(\xi\otimes\psi) 
 \;=\; 
 \big(d \,+\, \o_\M\big)(\xi\otimes\psi)
 \;=\; 0\,.
 $$
 Hence, $\xi\otimes\psi$ defines a parallel spinor on $\M$.
 \leer

 We remark that the parallel spinor $\xi\otimes\psi$ also gives rise
 to a Killing spinor on $S$. This follows of course from the general
 equivalence between Killing spinors on a Riemannian manifold and
 parallel spinors on its cone, as described by C. B\"ar (cf.~\cite{baer})
 and briefly summarized in appendix \ref{spinorsoncones}. Nevertheless,
 we will include this construction since it is an immediate consequence
 of our approach.
 We have seen above that $\xi\otimes\psi$, considered as function on 
 $P_\M$, satisfies $(d \,+\, \o_\M)(\xi\otimes\psi)
 \,=\; 0\, $. If we substitute  the expression $\o_\M$ given
 in Lemma \ref{omh} we obtain
 $$
 \big(d \,+\, \o_S \, + \, \sqrt{\tfrac\k{16n(n+2)}} \,
  \theta_S \, \wedge \, \1
 \big) (\xi\otimes\psi) \;=\; 0\,.
 $$
 Interpreted as equation on $P$  this is just the Killing equation, i.e.
 $$
 (d \, + \, \o_S) \,(\xi\otimes\psi) \, = \,
 -\,\tfrac{1}{2} \, \sqrt{\tfrac\k{16n(n+2)}}\,
    \theta_S \mathop{\cdot}\limits_S (\xi\otimes\psi) \, .
 $$
 The additional factor $1 \over 2$ is due to the action of the orthogonal
 Lie algebra on the spinor bundle.
 Before closing this section we formulate the lemma needed above.
 \begin{Lemma}\label{thetastar}
 $$
  \sqrt{\tfrac{\k}{8n(n+2)}}\,  
  \id\otimes\theta_{M}^{H\otimes E}\star
  \;=\; 
  \o_\M^{\H^n} \,.
 $$
 \end{Lemma}
 The proof needs an additional proposition which
 is analogous to Proposition 2.3 in \cite{qklast}.
 \begin{Proposition}
  Let $p$,\ $q$ be a base of $\C^2\cong\H$ with $\sigma(\,p,\,q\,)=1$ and 
  $f_1,f_2 \in F$. Then we have the following identity of operators on
  the spinor bundle:
  $$
   (p\otimes f_1)\wedge (q\otimes f_2) - (q\otimes f_1)\wedge (p\otimes f_2)
   \; = \; \id \otimes f_1\cdot f_2 \,.
  $$
  As element of $\so(4n+4)$, the left hand side acts on the spinor module
  via the isomorphism $\so(4n+4) \cong \spin(4n+4)$ sending $e_1\wedge e_2$
  to $\frac{1}{2}(\, e_1\;e_2 \,+\, \la\;e_1,\;e_2\;\ra\,)$.
 \end{Proposition}
 We remark that the proof of these two technical Propositions amounts to prove
 the decomposition (\ref{spinor}) of the spinor bundle.
 With the help of this proposition it is easy to prove Lemma \ref{thetastar}.
 \leer
 \proof
 Let $p$,\ $q$ be the canonical base of $\C^2$ used to define Clifford
 multiplication. We can then extend the isomorphism $\C\otimes_\R F\cong
 \C^2\otimes F$ in equation (\ref{FtoF}) to $\C\otimes_\R\L^2 F\cong
 \L^2(\C^2\otimes F)$ to get
 $$
 \begin{array}{rcl}
  (\theta_M \wedge \1) &=&
  \phantom{-}\frac{1}{2}\,(p\otimes \theta_M + q\otimes J\theta_M)\wedge
    (p\otimes \1 + q\otimes \J)\\[1.4ex]
  (i\theta_M \wedge \I) &=&
  -\frac{1}{2}\,(p\otimes \theta_M - q\otimes J\theta_M)\wedge
    (p\otimes \1 - q\otimes \J)\\[1.4ex]
  (j\theta_M \wedge \J) &=&
  \phantom{-}\frac{1}{2}\,(p\otimes J\theta_M - q\otimes \theta_M)\wedge
    (p\otimes \J - q\otimes \1)\\[1.4ex]
  (k\theta_M \wedge \K) &=&
  -\frac{1}{2}\,(p\otimes J\theta_M + q\otimes \theta_M)\wedge
    (p\otimes \J + q\otimes \1)\,.
 \end{array}
 $$
 In the second and fourth line we use that $J$ is conjugate linear. By
 summing up the four equations, some terms cancel, and we obtain
 $$
 \begin{array}{rcl}
  \o_\M^{\H^n} & = & \sqrt{\frac{\k}{16n(n+2)}}
  \,\Big((p\otimes \theta_M)\wedge(q\otimes \J) 
     - (q\otimes \theta_M)\wedge(p\otimes \J) \\[1.4ex]
   && \phantom{\sqrt{\frac{\k}{16n(n+2)}}
  \,\Big(}
    -(p\otimes J\theta_M)\wedge(q\otimes \1) 
    + (q\otimes J\theta_M)\wedge(p\otimes \1)\Big) \,.
 \end{array}
 $$
 Applying the proposition above yields the following equivalence
 $$
 \begin{array}{rcl}
  \o_\M^{\H^n} & = & \sqrt{\frac{\k}{16n(n+2)}}\, 
    (\id \otimes \theta_M\cdot \J - \id \otimes J\theta_M\cdot \1) \\[1.6ex]
  & = & \sqrt{\frac{\k}{16n(n+2)}}\, 
     \id \otimes (\J\otimes \theta_M - \1\otimes J\theta_M)\star
  \quad = \quad \sqrt{\frac{\k}{8n(n+2)}}\, 
     \id \otimes \theta_M^{H\otimes E}\star\,. \qed
 \end{array}
 $$

\section{Proof of the Theorem}

 The fact that a quaternionic Killing spinor on $M$ translates into a parallel
 spinor on the hyperk\"ahler manifold $\M$ is crucial to the proof of
 the main theorem.
 \setcounter{Theorem}{0}
 \begin{Theorem}
  Let $M$ be a compact quaternionic K\"ahler manifold of quaternionic
  dimension $n$ and positive scalar curvature $\k >0$. If there is an
  eigenspinor for the Dirac operator with eigenvalue $\l$ satisfying
  $$
   \l^2 = \frac{\k}{4}\frac{n+3}{n+2}\,,
  $$
  then $M$ is isometric to the quaternionic projective space.
 \end{Theorem}
 \leer
 \proof
 After the work done in the preceding sections the proof of this theorem
 reduces to a simple holonomy argument. The spinor bundle of the hyperk\"ahler
 manifold $\M$ is associated to the $\Sp(n)$--representation
 $$
  \S = \bigoplus_{r=0}^{n+1} \S_r
         = \bigoplus_{r=0}^{n+1} \Sym^r\C^2 \otimes \L_\circ^{n+1-r}F \,,
 $$
 where $F=H\oplus E$, and $\Sp(n)$ operates trivially on $H$.
 The holonomy of $\M$ is contained in $\Sp(n+1)$ and operates
 trivially on $\Sym^r\C^2$. The subbundle $\spb_{n+1}(\M)$ of $\spb(\M)$
 associated to $\S_{n+1}=\Sym^{n+1}\C^2$ is consequently trivialized by
 $n+2=\dim\,(\Sym^{n+1}\C^2)$ linearly independent parallel spinors.

 If $M$ admits a quaternionic Killing spinor, there are additional
 parallel spinors on $\M$, which are sections of the subbundle
 $\spb_1(\M)$ associated to the representation
 $\C^2\otimes \L_\circ^n F$. 

 Due to a result of Wang \cite{wang}, on a manifold with holonomy  equal to
 $\Sp(n+1)$ there are exactly $n+2$ linearly independent parallel spinors,
 just those trivializing $\spb_{n+1}(\M)$. The additional parallel
 spinor constructed out of a quaternionic Killing spinor reduces the holonomy
 further. According to Berger's list this can only happen if $\M$ is
 reducible or locally symmetric. In the first case, as consequence of a
 theorem of Gallot \cite{gallot}, $\M$ has to be flat. But so it is
 in the second case, because it is hyperk\"ahler, hence Ricci--flat, and in
 addition locally symmetric. Therefore $\M$ is flat which forces $M$
 to be isometric to the quaternionic projective space.
 \qed

\appendix
\section{Spinors on Cones}\label{spinorsoncones}

 In this appendix we will describe how to lift spinors on a Riemannian spin
 manifold $N$ to spinors on its cone $\widehat N := \R^+\times N$. In
 particular, we will show that Killing spinors on $N$ translate into parallel
 spinors on $\widehat N$. This construction is originally due to C.~B\"ar
 \cite{baer}.

 Let $(N,g_N)$ be a spin manifold of dimension $n$ and let $\pi: \widehat N =
 \R^+\times N \to N$ be the cone over $N$ endowed with the warped product
 metric $g_{\widehat N} = dt^2 + e^{2\l} \pi^* g_N$. The soldering and
 connection form on the principal bundle $P := P_{\Spin(n)}N$ associated
 to the chosen spin structure are denoted by $\theta_N$ and $\o_N$. Obviously,
 the cone is again spin and the spin structure reduces to the principal
 $\Spin(n)$--bundle $\widehat P:= \R^+\times P$. Forms on $P$ give rise
 to forms on $\widehat P$ by extending them constantly along the
 $\R^+$--direction. It easy to see that the soldering resp.~the
 connection form of $\widehat N$ on $\widehat P$ are given by
 $$
  \theta_{\widehat N} = dt + e^\l \theta_N \quad \mathrm{resp.}\quad
  \o_{\widehat N} = \o_N - e^\l \frac{\partial \l}{\partial t}\,
    \theta_N \wedge \Xi,
 $$
 where $\Xi = \frac{\partial}{\partial t}$ denotes the vertical unit vector.

 A spinor $\psi$ on $N$ can be interpreted as a $\Spin(n)$--equivariant
 function on $P$ with values in the spinor module $\S_n$ with a fixed
 Clifford module structure. With the canonical isomorphism
 $\CCl_n \cong \CCl_{n+1}^0$, $e_i \mapsto e_i\cdot \Xi$ in mind, we can
 consider $\S_n$ as an $\CCl_{n+1}^0$-- and therefore as an
 $\Spin(n+1)$--representation. The values of $\psi$ have now to be
 interpreted as lying in the $\Spin(n+1)$--representation. Let $\psi$
 be a Killing spinor on $N$. Interpreted as function on $P$, it satisfies:
 $$
  d\psi + (\o_N - \mu \theta_N)\, \psi \;=\; 0 \,,
 $$
 where $\mu$ is the Killing constant. Extending $\psi$ constantly in
 $\R^+$--direction and using the isomorphism of the Clifford algebras
 above defines a spinor on $\widehat{N}$. As a function on $\widehat P$,
 it satisfies
 $$
  d \psi + (\o_N - \mu \pi^*\theta_N \Xi)\, \psi \;=\; 0 \,.
 $$
 If the warping function is chosen such that
 $\mu = \frac{1}{2}e^\l \frac{\partial \l}{\partial t}$, the expression
 in brackets is
 equal to $ \o_{\widehat N} $, and therefore $\psi$ has to be parallel on
 $\widehat N$. The other way round, it is also clear that, if $\psi$ is a
 parallel spinor on $\widehat N$, then its associated function on
 $\widehat{P}$ is constant in $\R^+$--direction and it projects onto
 a Killing spinor on $N$.

\section{The Curvature Tensor of $\M$}\label{proofrhyper}

 In this appendix we will prove of Proposition \ref{hypcurv}
 relating the curvature tensor of $\M$ to the hyperk\"ahler
 part $R^{hyper}$ of the curvature of $M$. We recall that the connection
 form of $\M$ restricted to $P_\M$ is pulled back
 from $P$ and so is its curvature $\O_\M$. We have seen in
 section \ref{sixone} that the most convenient way to read the connection
 form $\o_\M$ given in Lemma \ref{omh} is
 $$
  \o_\M\,=\,\o_\M^{\sp(1)}
   \,+\,\o_\M^{\sp(n)}
   \,+\,\o_\M^{\H^n}\,,
 $$
 where $\o_\M^{\sp(1)}$ and $\o_\M^{\sp(n)}$ are the
 pull--backs of $\o_M^{\sp(1)}$ and $\o_M^{\sp(n)}$ respectively, and
 $$
  \o_\M^{\H^n} := \sqrt{\tfrac\k{16n(n+2)}}\;
  \big( \theta_M\wedge\1 + i\theta_M\wedge\I
   + j\theta_M\wedge\J + k\theta_M\wedge\K \big)\,.
 $$
 Defining $[\alpha\wedge\beta](\,X,\,Y\,):=[\alpha(X),\beta(Y)]-
 [\alpha(Y),\beta(X)]=[\beta\wedge\alpha](\,X,\,Y\,)$ for Lie algebra
 valued 1--forms $\alpha$,\ $\beta$, the curvature 2--form $\O_M$ of
 $M$ on $P$ can be written
 $$
 \begin{array}{rcl}
  \O_M
  & = & d\o_M + \tfrac{1}{2}\big[\o_M\wedge\o_M\big] \\[1.4ex]
  & = & \Big( d\o_M^{\sp(1)} + \tfrac{1}{2}\big[\o_M^{\sp(1)}
   \wedge\o_M^{\sp(1)}\big] \Big) + \Big( d\o_M^{\sp(n)} +
   \tfrac{1}{2}\big[\o_M^{\sp(n)}\wedge\o_M^{\sp(n)}\big] \Big)\,.
 \end{array}
 $$
 Of course, $[\o_M^{\sp(1)}\wedge\o_M^{\sp(n)}]=0$ since $\sp(1)$
 and $\sp(n)$ centralize each other in $\sp(1)\oplus\sp(n)$.
 Using the naturality of the exterior differential and
 $\o_M^{\sp(1)}\,=\,\pi_S^*\o$ we conclude that the first
 summand is equal to $\pi_S^*\O$. Thus, it corresponds to
 $-\tfrac\k{8n(n+2)}\;R^H$ in the sense of decomposition
 (\ref{curv1}) as shown in the proof of Lemma \ref{curv}. We conclude
 that the second summand corresponds to $-\tfrac\k{8n(n+2)}\;R^E
 \;+\;R^{hyper}$.
 The curvature $\O_\M$ of $\M$ can be calculated
 similarly. With $\sp(1)$ and $\sp(n)$ centralizing each other in
 $\sp(n+1)$ we still have $[\o_\M^{\sp(1)}\wedge\o_\M^{\sp(n)}]=0$,
 and $\O_\M$ is the sum of the three terms:
 \be\label{curvmh}
 \begin{array}{ccccc}
  d\o_\M^{\sp(1)} &\!\!+\!\!&
  \tfrac{1}{2}\;\big[\o_\M^{\sp(1)}\wedge\o_\M^{\sp(1)}\big] &\!\!+\!\!&
  \tfrac{1}{2}\;\big[\o_\M^{\H^n}\wedge\o_\M^{\H^n}\big]^{\sp(1)}
  \\[1ex]
  d\o_\M^{\H^n} &\!\!+\!\!&
  \,\,\,\,\,\big[\o_\M^{\sp(1)}\wedge\o_\M^{\H^n}\big] &\!\!+\!\!&
  \!\!\!\!\!\big[\o_\M^{\sp(n)}\wedge\o_\M^{\H^n}\big]
  \\[1ex]
  d\o_\M^{\sp(n)} &\!\!+\!\!&
  \tfrac{1}{2}\;\big[\o_\M^{\sp(n)}\wedge\o_\M^{\sp(n)}\big] &\!\!+\!\!&
  \tfrac{1}{2}\;\big[\o_\M^{\H^n}\wedge\o_\M^{\H^n}\big]^{\sp(n)}
 \end{array}
 \ee
 where $[\o_\M^{\H^n}\wedge\o_\M^{\H^n}]$ is projected onto its two
 components in $\sp(1)$ and $\sp(n)$ according to the Cartan
 decomposition $\sp(n+1)=(\sp(1)\oplus\sp(n))\oplus\H^n$.
 Using the formula
 $$
  [\,a_1\wedge b_1,\, a_2\wedge b_2\,]
  = \la a_1,a_2\ra b_1\wedge b_2 - \la b_1,a_2\ra a_1\wedge b_2
  - \la a_1,b_2\ra b_1\wedge a_2 + \la b_1,b_2\ra a_1\wedge a_2
 $$
 we find
 \begin{eqnarray}
  \lefteqn{\tfrac{1}{2}\big[\o_\M^{\H^n}\wedge
   \o_\M^{\H^n}\big]} && \nonumber\\
  &\!\!= &
  \tfrac\k{32n(n+2)}\Big(
   \theta_M\wedge \theta_M \,+\, i\theta_M\wedge i\theta_M \,+\,
   j\theta_M\wedge j\theta_M \,+\, k\theta_M\wedge k\theta_M \nonumber\\
  & & \qquad\qquad
  \,+\, 2\la\theta_M\wedge i\theta_M\ra \1\wedge\I
  \,+\, 2\la\theta_M\wedge j\theta_M\ra \1\wedge\J
  \,+\, 2\la\theta_M\wedge k\theta_M\ra \1\wedge\K \nonumber\\
  & & \qquad\qquad
  \,+\, 2\la i\theta_M\wedge j\theta_M\ra \I\wedge\J
  \,+\, 2\la i\theta_M\wedge k\theta_M\ra \I\wedge\K
  \,+\, 2\la j\theta_M\wedge k\theta_M\ra \J\wedge\K
  \Big)\nonumber\\
  &\!\!= & \label{RE}
  \tfrac\k{32n(n+2)}\Big(
   \theta_M\wedge \theta_M \,+\, i\theta_M\wedge i\theta_M \,+\,
   j\theta_M\wedge j\theta_M \,+\, k\theta_M\wedge k\theta_M
  \Big) \\
  & & \label{RH} + 
  \tfrac\k{16n(n+2)}\Big(
   \la\theta_M\!\wedge\! i\theta_M\ra(\1\!\wedge\!\I\!-\!\J\!\wedge\!\K) +
   \la\theta_M\!\wedge\! j\theta_M\ra(\1\!\wedge\!\J\!-\!\K\!\wedge\!\I) +
   \la\theta_M\!\wedge\! k\theta_M\ra(\1\!\wedge\!\K\!-\!\I\!\wedge\!\J)
  \Big)\qquad
 \end{eqnarray}
 We remark that by construction of the base $\1$,\ $\I$,\ $\J$ and $\K$
 the infinitesimal action of $i$,\ $j$ and $k\in\sp(1)$ on $H$ can be written 
 $$
  \begin{array}{rlclcrlclcrlcl}
   i:\!\! & \1\mapsto -\I & \!\!\! & \J\mapsto\phantom{-}\K & \quad
   j:\!\! & \1\mapsto -\J & \!\!\! & \K\mapsto\phantom{-}\I & \quad
   k:\!\! & \1\mapsto -\K & \!\!\! & \I\mapsto\phantom{-}\J \\
   & \I\mapsto\phantom{-}\1 & & \K\mapsto -\J &
   & \J\mapsto\phantom{-}\1 & & \I\mapsto -\K &
   & \K\mapsto\phantom{-}\1 & & \J\mapsto -\I
  \end{array}\,,
 $$
 i.~e.~$i$ corresponds to $-(\1\wedge\I -\J\wedge\K)$. Thus, the
 summand (\ref{RH}) is equal to $-\pi_S^*\O$ according to Lemma \ref{curv},
 i.~e.~equal to $\tfrac\k{8n(n+2)}R^H$. Without proof we state that the
 summand (\ref{RE}) is equal to $\tfrac\k{8n(n+2)}R^E$ as this is
 certainly true on the quaternionic projective space. Consequently,
 in the decomposition (\ref{curvmh}) of $\O_\M$ the first
 term vanishes as does the second, because straightforward
 calculations show that it depends linearly on the torsion of the
 Levi--Civit\'a connection $\o_M$ of $M$. Hence, the curvature
 $\O_\M$ reduces to 
 $$
  \O_\M = d\o_\M^{\sp(n)} \,+\,
  \tfrac{1}{2}\;\big[\o_\M^{\sp(n)}\wedge\o_\M^{\sp(n)}\big] \,+\,
  \tfrac\k{8n(n+2)} R^E = R^{hyper}\,.
 $$
 In this way $\O_\M$ operates only on the subbundle
 $P_\M\times_{\Sp(n)}E$ of the tangent bundle of $\M$.
 We remark that this subbundle is canonically isomorphic to $P_\M
 \times_{\Sp(n)}T$, i.~e.~to $\widehat{\pi}^*TM$. \qed

\end{document}